\documentclass{amsart}

\usepackage[a4paper,hmargin=3.5cm,vmargin=4cm]{geometry}
\usepackage{amsfonts,amssymb,amscd,amstext}
\usepackage{graphicx}
\usepackage[dvips]{epsfig}
\usepackage{mathpazo}

\renewcommand{\leq}{\leqslant}
\renewcommand{\geq}{\geqslant}
\newcommand{\ptl}{\partial}
\newcommand{\rr}{{\mathbb{R}}}
\newcommand{\la}{\lambda}
\newcommand{\hh}{{\mathbb{H}}}
\newcommand{\sph}{{\mathbb{S}}}
\newcommand{\stres}{{\mathbb{S}^3}}
\newcommand{\h}{\mathcal{H}}
\newcommand{\sub}{\subset}
\newcommand{\subeq}{\subseteq}
\newcommand{\escpr}[1]{\big<#1\big>}
\newcommand{\Sg}{\Sigma}
\newcommand{\sg}{\sigma}
\newcommand{\Om}{\Omega}
\newcommand{\eps}{\varepsilon}
\newcommand{\var}{\varphi}
\newcommand{\ga}{\gamma}
\newcommand{\Ga}{\Gamma}

\DeclareMathOperator{\divv}{div}

\newtheorem{theorem}{Theorem}[section]
\newtheorem{proposition}[theorem]{Proposition}
\newtheorem{lemma}[theorem]{Lemma}

\theoremstyle{definition}

\newtheorem{remark}[theorem]{Remark}
\newtheorem{example}[theorem]{Example}

\theoremstyle{remark}

\numberwithin{equation}{section}

\begin{document}

\title[Area-stationary surfaces in the sub-Riemannian
$\sph^3$]{area-stationary surfaces inside \\ the sub-riemannian
three-sphere}

\author[A.~Hurtado]{Ana Hurtado}
\address{Departament de Matem\`atiques \\ Universitat Jaume I
\\ 8029 AP Castell\'o \\ Spain}
\email{ahurtado@mat.uji.es}

\author[C.~Rosales]{C\'esar Rosales}
\address{Departamento de Geometr\'{\i}a y Topolog\'{\i}a \\
Universidad de Granada \\ E--18071 Granada \\ Spain}
\email{crosales@ugr.es}

\date{August 1, 2006}

\thanks{The first author has been partially supported by MCyT-Feder
research project MTM2004-06015-C02-01.  The second author has
been supported by MCyT-Feder research project MTM2004-01387}

\begin{abstract}
We consider the sub-Riemannian metric $g_{h}$ on $\stres$ provided by the restriction of the
Riemannian metric of curvature $1$ to the plane distribution orthogonal to the Hopf vector
field.  We compute the geodesics associated to the Carnot-Carath\'eodory distance and we show
that, depending on their curvature, they are closed or dense subsets of a Clifford torus.

We study area-stationary surfaces with or without a volume constraint in $(\stres,g_{h})$.  By
following the ideas and techniques in \cite{stationary} we introduce a variational notion of
mean curvature, characterize stationary surfaces, and prove classification results for complete
volume-preserving area-stationary surfaces with non-empty singular set.  We also use the
behaviour of the Carnot-Carath\'eodory geodesics and the ruling property of constant mean
curvature surfaces to show that the only $C^2$ compact, connected, embedded surfaces in
$(\stres,g_{h})$ with empty singular set and constant mean curvature $H$ such that
$H/\sqrt{1+H^2}$ is an irrational number, are Clifford tori.  Finally we describe which are the
complete rotationally invariant surfaces with constant mean curvature in $(\stres,g_{h})$.
\end{abstract}

\maketitle

\thispagestyle{empty}

\section{Introduction}
\label{sec:intro}
\setcounter{equation}{0}

Sub-Riemannian geometry studies spaces equipped with a path metric
structure where motion is only possible along certain trajectories
known as \emph{admissible} (or \emph{horizontal}) curves.  This
discipline has motivations and ramifications in several parts of
mathematics and phy\-sics, such as Riemannian and contact geometry,
control theory, and classical mechanics.

In the last years the interest in variational questions in sub-Riemannian geometry has
increased.  One of the reasons for the recent growth of this field has been the desire to solve
global problems involving the sub-Riemannian area in the Heisenberg group. The $3$-dimensional
Heisenberg group $\hh^1$ is one of the simplest and most important non-trivial sub-Riemannian
manifolds, and it is object of an intensive study.  In fact, some of the classical
area-minimizing questions in Euclidean space such as the Plateau problem, the Bernstein
problem, or the isoperimetric problem have been treated in $\hh^1$.  Though these problems are
not completely solved, some important results have been established, see \cite{pauls},
\cite{chy}, \cite{chmy}, \cite{stationary}, \cite{survey}, and the references therein.  For
example in \cite{stationary}, M.~Ritor\'e and the second author have proved that the only $C^2$
isoperimetric solutions in $\hh^1$ are the spherical sets conjectured by P.~Pansu \cite{pansu}
in the early eighties.  The particular case of $\hh^1$ has inspired the study of similar
questions as that as the development of a theory of constant mean curvature surfaces in
different classes of sub-Riemannian manifolds, such as Carnot groups \cite{dgn}, see also
\cite{dgn2}, pseudohermitian manifolds \cite{chmy}, vertically rigid manifolds \cite{hp}, and
contact manifolds \cite{pre}.

Besides the Heisenberg group, one of the most important examples in sub-Riemannian geometry
comes from the Heisenberg spherical structure, see \cite{Gr} and \cite[\S~11]{montgomery}. In
this paper we use the techniques and arguments employed in \cite{stationary} to study
\emph{area-stationary surfaces with or without a volume constraint inside the sub-Riemannian
$3$-sphere}. Let us precise the situation.  We denote by $(\stres,g)$ the unit $3$-sphere
endowed with the Riemannian metric of constant sectional curvature $1$. This manifold is a
compact Lie group when we consider the quaternion product $p\cdot q$.  A basis of right
invariant vector fields in $(\stres,\cdot)$ is given by $\{E_1,E_2,V\}$, where $E_1(p)=j\cdot
p$, $E_2(p)=k\cdot p$ and $V(p)=i\cdot p$ (here $i$, $j$ and $k$ are the complex quaternion
units).  The vector field $V$ is sometimes known as the \emph{Hopf vector field} in $\stres$
since its integral curves parameterize the fibers of the Hopf map
$\mathcal{F}:\stres\to\sph^2$.  We equip $\stres$ with the \emph{sub-Riemannian metric} $g_{h}$
provided by the restriction of $g$ to the \emph{horizontal distribution}, which is the smooth
plane distribution generated by $E_1$ and $E_2$.  Inside the sub-Riemannian manifold
$(\stres,g_h)$ we can consider many of the notions existing in Riemannian geometry. In
particular, we can define the \emph{Carnot-Carath\'eodory distance} $d(p,q)$ between two
points, the \emph{volume} $V(\Om)$ of a Borel set $\Om$, and the \emph{area} $A(\Sg)$ of a
$C^1$ immersed surface $\Sg$, see Section~\ref{sec:preliminaries} and the beginning of
Section~\ref{sec:geodesics} for precise definitions.

In Section~\ref{sec:geodesics} we use intrinsic arguments similar to those in \cite[\S
3]{stationary} to study \emph{geodesics} in $(\stres,g_h)$.  They are defined as $C^2$
horizontal curves which are critical points of the Riemannian length for variations by
horizontal curves with fixed extreme points.  Here ``horizontal" means that the tangent vector
to the curve lies in the horizontal distribution.  The geodesics are solutions of a second
order linear differential equation depending on a real parameter called the $\emph{curvature}$
of the geodesic, see Proposition~\ref{prop:eqgeo}.  As was already observed in \cite{chmy} the
geodesics of curvature zero coincide with the horizontal great circles of $\stres$.  From an
explicit expression of the geodesics we can easily see that they are horizontal lifts via the
Hopf map $\mathcal{F}:\stres\to\sph^2$ of the circles of revolution in $\sph^2$. Moreover, in
Proposition~\ref{prop:geoclosed} we show that the topological behaviour of a geodesic $\ga$
only depends on its curvature $\la$. Precisely, if $\la/\sqrt{1+\la^2}$ is a rational number
then $\ga$ is a closed curve diffeomorphic to a circle.  Otherwise $\ga$ is diffeomorphic to a
straight line and coincides with a dense subset of a Clifford torus in $\stres$.  We finish
Section~\ref{sec:geodesics} with the notion of $\emph{Jacobi field}$ in $(\stres,g_h)$.  These
vector fields are associated to a variation of a given geodesic by geodesics of the same
curvature.  They will be key ingredients in some proofs of Section~\ref{sec:classification}.

In Section~\ref{sec:mainresult} we consider critical surfaces with or without a volume
constraint for the area functional in $(\stres,g_h)$. These surfaces have been well studied in
the Heisenberg group $\hh^1$, and most of their properties remain valid, with minor
modifications, in the sub-Riemannian $3$-sphere.  For example if $\Sg$ is a $C^2$
volume-preserving area-stationary surface then the \emph{mean curvature} of $\Sg$ defined in
\eqref{eq:mc} is constant off of the \emph{singular set} $\Sg_0$, the set of points where the
surface is tangent to the horizontal distribution. Moreover $\Sg-\Sg_0$ is a ruled surface in
$(\stres,g_h)$ since it is foliated by geodesics of the same curvature.  Furthermore, by the
results in \cite{chmy}, the singular set $\Sg_0$ consists of isolated points or $C^1$ curves.
We can also prove a characterization theorem similar to \cite[Thm.~4.16]{stationary}: for a
$C^2$ surface $\Sg$, to be area-stationary with or without a volume constraint is equivalent to
that $H$ is constant on $\Sg-\Sg_0$ and the geodesics contained in $\Sg-\Sg_0$ meet
\emph{orthogonally} the singular curves.  Though the proofs of these results are the same as in
\cite{stationary} we state them explicitly since they are the starting points to prove our
classification results in Section~\ref{sec:classification}.

In \cite{chmy}, J.-H.~Cheng, J.-F. Hwang, A.~Malchiodi and P.~Yang
found the first examples of constant mean curvature surfaces in
$(\stres,g_h)$.  They are the totally geodesic $2$-spheres in
$(\stres,g)$ and the Clifford tori $\mathcal{T}_\rho$ defined in
complex notation by the points $(z_1,z_2)\in\stres$ such that
$|z_1|^2=\rho^2$.  The above mentioned authors also established two
interesting results for compact surfaces with constant mean curvature
in $(\stres,g_h)$.  First they gave a strong topological restriction
by showing \cite[Thm.  E]{chmy} that such a surface must be
homeomorphic either to a sphere or to a torus.  Second they obtained
\cite[Proof of Cor.  F]{chmy} that any compact, embedded, $C^2$
surface with vanishing mean curvature and at least one isolated
singular point must coincide with a totally geodesic $2$-sphere in
$(\stres,g)$.

\begin{figure}[h]
\centering{\includegraphics[width=5.5cm]{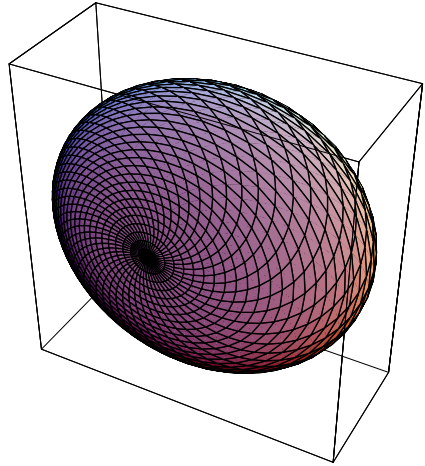}}
\caption{Stereographic projection from $\stres-\{\text{north pole}\}$
to $\rr^3$ of a spherical surface $\sph_{\la}$ given by the union of
all the geodesics of curvature $\la$ and length $\pi/\sqrt{1+\la^2}$
leaving from $p=(1,0,0,0)$.}
\label{fig:spheres}
\end{figure}

In Section~\ref{sec:classification} of the paper we give the complete
classification of complete, volume-preser\-ving area-stationary
surfaces in $(\stres,g_h)$ with non-empty singular set.  In
Theorem~\ref{th:spheres} we generalize the aforementioned Theorem E in
\cite{chmy}: we prove that if $\Sg$ is a $C^2$ complete, connected,
immersed surface with constant mean curvature $H$ and at least an
isolated singular point, then $\Sg$ is congruent with the spherical
surface $\mathcal{S}_H$ described as the union of all the geodesics of
curvature $H$ and length $\pi/\sqrt{1+H^2}$ leaving from a given
point, see Figure~\ref{fig:spheres}.  Our main result in this section
characterizes complete volume-preserving area-stationary surfaces with
at least one singular curve $\Ga$.  The local description given in
Theorem~\ref{th:chmy} of such a surface $\Sg$ around $\Ga$, and the
orthogonality condition between singular curves and geodesics in
Theorem~\ref{th:constant}, imply that a small neighborhood of $\Ga$ in
$\Sg$ consists of the union of small pieces of all the geodesics
$\ga_\eps$ of the same curvature leaving from $\Ga$ orthogonally.  By
using the completeness of $\Sg$ we can extend these geodesics until
they meet another singular point.  Finally, from a detailed study of
the Jacobi vector field associated to the family $\ga_\eps$ we deduce
that the singular curve $\Ga$ must be a geodesic in $(\stres,g_h)$.
This allows us to conclude that $\Sg$ is congruent with one of the
surfaces $\mathcal{C}_{\mu,\la}$ obtained when we leave orthogonally
from a given geodesic of curvature $\mu$ by geodesics of curvature
$\la$, see Example~\ref{ex:smula}.

The classification of complete surfaces with empty singular set and
constant mean curvature in $(\stres,g_h)$ seems to be a difficult
problem.  In Section~\ref{sec:classification} we prove some results in
this direction.  In Proposition~\ref{prop:vertical} we show that the
Clifford tori $\mathcal{T}_\rho$ are the only complete surfaces with
constant mean curvature such that the Hopf vector field $V$ is always
tangent to the surface.  In Theorem~\ref{th:irrat} we characterize the
Clifford tori as the unique compact embedded surfaces with empty
singular set and constant mean curvature $H$ such that
$H/\sqrt{1+H^2}$ is irrational.  These results might suggest that
Theorem~\ref{th:irrat} holds without any further assumption on the
curvature $H$ of the surface.

In the last section of the paper we describe complete surfaces with
constant mean curvature in $(\stres,g_h)$ which are invariant under
the isometries of $(\stres,g)$ fixing the vertical equator passing
through $(1,0,0,0)$.  For such a surface the equation of constant mean
curvature can be reduced to a system of ordinary differential
equations.  Then, a detailed analysis of the solutions yields a
counterpart in $(\stres,g_h)$ of the classification by C.~Delaunay of
rotationally invariant constant mean curvature surfaces in $\rr^3$,
later extended by W.-H.~Hsiang \cite{delaunay} to $(\stres,g)$.  In
particular we can find compact, embedded, unduloidal type surfaces
with empty singular set and constant mean curvature $H$ such that
$H/\sqrt{1+H^2}$ is rational.  This provides an example illustrating
that all the hypotheses in Theorem~\ref{th:irrat} are necessary.

In addition to the geometric interest of this work we believe that our results may be applied
in two directions.  First, they could be useful to solve the \emph{isoperimetric problem} in
$(\stres,g_h)$ which consists of enclosing a fixed amount of volume with the least possible
boundary area.  In fact, if we assume that the solutions to this problem are $C^2$ smooth and
have at least one singular point, then they must coincide with one of the surfaces
$\mathcal{S}_\la$ or $\mathcal{C}_{\mu,\la}$ introduced in Section~\ref{sec:classification}.
Second, our classification results could be utilized to find examples of constant mean
curvature surfaces inside the Riemannian Berger spheres $(\stres,g_k)$.  This is motivated by
the fact that the metric space $(\stres,d)$ associated to the Carnot-Carath\'eodory distance is
limit, in the Gromov-Hausdorff sense, of the spaces $(\stres,d_k)$, where $d_k$ is the
Riemannian distance of $g_k$ \cite[p.~109]{Gr}.

The authors want to express their gratitude to O.~Gil and M.~Ritor\'e
for encouraging them to write these notes and helping discussions.
This work was initiated while A.~Hurtado was visiting the University
of Granada in the winter of 2006.  The paper was finished during a
short visit of C.~Rosales to the University Jaume I (Castell\'o) in
the summer of 2006.

\section{Preliminaries}
\label{sec:preliminaries}
\setcounter{equation}{0}

Throughout this paper we will identify a point
$p=(x_{1},y_{1},x_{2},y_{2})\in\rr^4$ with the quaternion
$x_{1}+iy_{1}+jx_{2}+ky_{2}$.  We denote the quaternion product and
the scalar product of $p,q\in\rr^4$ by $p\cdot q$ and $\escpr{p,q}$,
respectively.  The unit sphere $\stres\sub\rr^4$ endowed with the
quaternion product is a compact, noncommutative, $3$-dimensional Lie
group.  For $p\in\stres$, the \emph{right translation} by $p$ is the
diffeomorphism $R_p(q)=q\cdot p$.  A basis of right invariant vector
fields in $(\stres,\cdot)$ given in terms of the Euclidean coordinate
vector fields is
\begin{align*}
V(p):&=i\cdot p=-y_{1}\,\frac{\ptl}{\ptl x_{1}}+x_{1}\,\frac{\ptl}{\ptl y_{1}}
-y_{2}\,\frac{\ptl}{\ptl x_{2}}+x_{2}\,\frac{\ptl}{\ptl y_{2}},
\\
E_{1}(p):&=j\cdot p=-x_{2}\,\frac{\ptl}{\ptl x_{1}}+y_{2}\,\frac{\ptl}{\ptl
y_{1}}+x_{1}\,\frac{\ptl}{\ptl x_{2}}-y_{1}\,\frac{\ptl}{\ptl y_{2}},
\\
E_{2}(p):&=k\cdot p=-y_{2}\,\frac{\ptl}{\ptl x_{1}}-x_{2}\,\frac{\ptl}{\ptl
y_{1}}+y_{1}\,\frac{\ptl}{\ptl x_{2}}+x_{1}\,\frac{\ptl}{\ptl y_{2}}.
\end{align*}

We define the \emph{horizontal distribution} $\mathcal{H}$ in $\stres$
as the smooth plane distribution generated by $E_{1}$ and $E_{2}$.
The \emph{horizontal projection} of a vector $X$ onto $\mathcal{H}$ is
denoted by $X_{h}$.  A vector field $X$ is \emph{horizontal} if
$X=X_h$.  A \emph{horizontal curve} is a piecewise $C^1$ curve such
that the tangent vector (where defined) lies in the horizontal
distribution.

We denote by $[X,Y]$ the Lie bracket of two $C^1$ tangent vector
fields $X,Y$ on $\stres$.  Note that $[E_{1},V]=2E_{2}$,
$[E_{2},V]=-2E_{1}$ and $[E_{1},E_{2}]=-2V$, so that $\mathcal{H}$ is
a \emph{bracket generating distribution}.  Moreover, by Frobenius
theorem we have that $\mathcal{H}$ is nonintegrable.  The vector
fields $E_{1}$ and $E_{2}$ generate the kernel of the contact $1$-form
given by the restriction to the tangent bundle $T\stres$ of
$\omega:=-y_{1}\,dx_{1}+x_{1}\,dy_{1}-y_{2}\,dx_{2}+x_{2}\,dy_{2}$.

We introduce a \emph{sub-Riemannian metric} $g_h$ on $\stres$ by
considering the Riemannian metric on $\mathcal{H}$ such that
$\{E_{1},E_{2}\}$ is an orthonormal basis at every point.  It is
immediate that the Riemannian metric
$g=\escpr{\cdot\,,\cdot}|_{\stres}$ provides an extension to
$T\stres$ of the sub-Riemannian metric such that
$\{E_{1},E_{2},V\}$ is orthonormal.  The metric $g$ is
bi-invariant and so the right translations $R_{p}$ and the left
translations $L_{p}$ are isometries of $(\stres,g)$.  We denote by
$D$ the Levi-Civit\'a connection on $(\stres,g)$. The following
derivatives can be easily computed
\begin{alignat}{2}
\notag
D_{E_{1}}E_{1}&=0, \qquad \ \ \,\, D_{E_{2}}E_{2}=0,
\qquad \quad \ \,D_{V}V=0,
\\
\label{eq:christoffel}
D_{E_{1}}E_{2}&=-V, \qquad D_{E_{1}}V=E_{2}, \qquad \, \ \,\,
D_{E_{2}}V=-E_{1},
\\
\notag
D_{E_{2}}E_{1}&=V, \qquad \quad\!\!\,  D_{V}E_{1}=-E_{2}, \qquad \!
D_{V}E_{2} =E_{1} .
\end{alignat}

For any tangent vector field $X$ on $\stres$ we define $J(X):=D_XV$.
Then we have $J(E_{1})=E_{2}$, $J(E_{2})=-E_{1}$ and $J(V)=0$, so that
$J^2=-\text{Identity}$ when restricted to the horizontal distribution.
It is also clear that
\[
\escpr{J(X),Y}+\escpr{X,J(Y)}=0,
\]
for any pair of vector fields $X$ and $Y$.  The involution $J:\h\to\h$
together with the contact $1$-form
$\omega=-y_{1}\,dx_{1}+x_{1}\,dy_{1}-y_{2}\,dx_{2}+x_{2}\,dy_{2}$
provides a \emph{pseudohermitian structure} on $\stres$, as stated in
\cite[Appendix]{chmy}.  We remark that $J:\h\to\h$ coincides with the
restriction to $\h$ of the complex structure on $\rr^4$ given by the
left multiplication by $i$, that is
\[
J(X)=i\cdot X,\quad\text{ for any  } X\in\h.
\]

Now we introduce notions of volume and area in $(\stres,g_h)$.  We
will follow the same approach as in \cite{revolucion} and
\cite{stationary}.  The volume $V(\Om)$ of a Borel set
$\Om\subeq\stres$ is the Haar measure associated to the quaternion
product, which turns out to coincide with the Riemannian volume of
$g$.  Given a $C^1$ surface $\Sg$ immersed in $\stres$, and a unit
vector field $N$ normal to $\Sg$ in $(\stres,g)$, we define the area
of $\Sg$ in $(\stres,g_h)$ by
\begin{equation}
\label{eq:area}
A(\Sg):=\int_{\Sg}|N_{h}|\,d\Sg,
\end{equation}
where $N_{h}=N-\escpr{N,V}\,V$, and $d\Sg$ is the Riemannian area
element on $\Sg$.  If $\Om$ is an open set of $\stres$ bounded by a
$C^2$ surface $\Sg$ then, as a consequence of the Riemannian
divergence theorem, we have that $A(\Sg)$ coincides with the
sub-Riemannian perimeter of $\Om$ defined by
\[
\mathcal{P}(\Om)=\sup\,\left\{\int_{\Om}\divv X\,dv;\,|X|\leq 1\right\},
\]
where the supremum is taken over $C^1$ \emph{horizontal tangent
vector fields} on $\stres$. In the definition above $dv$ and
$\divv$ are the Riemannian volume and divergence of $g$,
respectively.

For a $C^1$ surface $\Sg\sub\stres$ the \emph{singular set} $\Sg_0$
consists of those points $p\in\Sg$ for which the tangent plane
$T_p\Sg$ coincides with $\h_{p}$.  As $\Sg_0$ is closed and has empty
interior in $\Sg$, the \emph{regular set} $\Sg-\Sg_0$ of $\Sg$ is open
and dense in $\Sg$.  It follows from the arguments in \cite[Lemme
1]{d2}, see also \cite[Theorem 1.2]{balogh}, that for a $C^2$ surface
$\Sg$ the Hausdorff dimension of $\Sg_{0}$ with respect to the
Riemannian distance in $\stres$ is less than two.  If $\Sg$ is
oriented and $N$ is a unit normal vector to $\Sg$ then we can describe
the singular set as $\Sg_{0}=\{p\in\Sg:N_h(p)=0\}$.  In the regular
part $\Sg-\Sg_0$, we can define the \emph{horizontal Gauss map}
$\nu_h$ and the \emph{characteristic vector field} $Z$, by
\begin{equation}
\label{eq:nuh}
\nu_h:=\frac{N_h}{|N_h|},\qquad Z:=J(\nu_h)=i\cdot\nu_{h}.
\end{equation}
As $Z$ is horizontal and orthogonal to $\nu_h$, we conclude that $Z$
is tangent to $\Sg$.  Hence $Z_{p}$ generates $T_{p}\Sg\cap\h_p$.  The
integral curves of $Z$ in $\Sg-\Sg_0$ will be called
\emph{characteristic curves} of $\Sg$.  They are both tangent to $\Sg$
and horizontal.  Note that these curves depend on the unit normal $N$
to $\Sg$.  If we define
\begin{equation}
\label{eq:ese}
S:=\escpr{N,V}\,\nu_h-|N_h|\,V,
\end{equation}
then $\{Z_{p},S_{p}\}$ is an orthonormal basis of $T_p\Sg$ whenever
$p\in\Sg-\Sg_0$.

Any isometry of $(\stres,g)$ leaving invariant the horizontal
distribution preserves the area $A(\Sg)$ of surfaces in
$(\stres,g_h)$.  Examples of such isometries are left and right
translations.  The rotation of angle $\theta$ given by
\begin{equation}
\label{eq:rtheta}
r_{\theta}(x_{1},y_{1},x_{2},y_{2})=(x_{1},y_{1},
(\cos\theta)x_{2}-(\sin\theta)y_{2},(\sin\theta)x_{2}+(\cos\theta)y_{2})
\end{equation}
is also such an isometry since it transforms the orthonormal basis
$\{E_{1},E_{2},V\}$ at $p$ into the orthonormal basis
$\{(\cos\theta)E_{1}
+(\sin\theta)E_{2},(-\sin\theta)E_{1}+(\cos\theta)E_{2},V\}$ at
$r_{\theta}(p)$.  We say that two surfaces $\Sg_{1}$ and $\Sg_{2}$ are
\emph{congruent} if there is an isometry $\phi$ of $(\stres,g)$
preserving the horizontal distribution and such that
$\phi(\Sg_{1})=\Sg_{2}$.

Finally we recall that the Hopf fibration
$\mathcal{F}:\stres\to\sph^2\equiv\stres\cap\{x_1=0\}$ is the Riemannian submersion given by
$\mathcal{F}(p)=\overline{p}\cdot i\cdot p$ (here $\overline{p}$ denotes the conjugate of the
quaternion $p$).  In terms of Euclidean coordinates we get
\[
\mathcal{F}(x_{1},y_{1},x_{2},y_{2})=(0,x^2_{1}+y^2_{1}-x_{2}^2-y_{2}^2,2\,
(x_{2}y_{1}-x_{1}y_{2}),2\,(x_{1}x_{2}+y_{1}y_{2})).
\]
The fiber passing through $p\in\stres$ is the great circle parameterized by $\exp(it)\cdot p$.
Clearly the fibers are integral curves of the vertical vector $V$, which is sometimes known as
the Hopf vector field.  A lift of a curve $c:(-\eps,\eps)\to\sph^2$ is a curve
$\ga:(-\eps,\eps)\to\stres$ such that $\mathcal{F}(\ga)=c$.  By general properties of principal
bundles we have that for any piecewise $C^1$ curve $c$ there is a unique horizontal lift of $c$
passing through a fixed point $p\in\mathcal{F}^{-1}(c(0))$, see \cite[p.~88]{kn}.  For any
$\rho\in (0,1)$ let $c_{\rho}$ be the geodesic circle of $\sph^2$ contained in the plane
$\{x_{1}=0,\,y_1=2\rho^2-1\}$.  The set $\mathcal{T}_{\rho}=\mathcal{F}^{-1}(c_{\rho})$ is the
Clifford torus in $\stres$ described by the pairs of complex numbers $(z_{1},z_{2})$ such that
$|z_{1}|^2=\rho^2$ and $|z_{2}|^2=1-\rho^2$.

\section{Carnot-Carath\'eodory geodesics in $\stres$}
\label{sec:geodesics}
\setcounter{equation}{0}

Let $\ga:I\to\stres$ be a piecewise $C^1$ curve defined on a compact interval $I\sub\rr$.  The
\emph{length} of $\ga$ is the Riemannian length $L(\ga):=\int_{I}|\dot{\ga}|$.  For any two
points $p,q\in\stres$ we can find, by Chow's connectivity theorem \cite[\S 1.2.B] {Gr}, a
$C^\infty$ horizontal curve joining these points.  The \emph{Carnot-Carath\'eodory distance}
$d(p,q)$ is defined as the infimum of the lengths of all piecewise $C^1$ horizontal curves
joining $p$ and $q$.  The topologies on $\stres$ defined by $d$ and the Riemannian distance
associated to $g$ are the same, see \cite[Cor. 2.6]{andre}.  In the metric space $(\stres,d)$
there is a natural extension for continuous curves of the notion of length, see
\cite[p.~19]{andre}.  We say that a continuous curve $\ga$ joining $p$ and $q$ is
\emph{length-minimizing} if $L(\ga)=d(p,q)$.  Since the metric space $(\stres,d)$ is complete
we can apply the Hopf-Rinow theorem in sub-Riemannian geometry \cite[Thm. 2.7]{andre} to ensure
the existence of length-minimizing curves joining two given points. Moreover, by \cite[Cor.
6.2]{geo}, see also \cite[Chapter~5]{montgomery}, any of these curves is $C^\infty$. In this
section we are interested in smooth curves which are critical points of length under any
variation by horizontal curves with fixed endpoints.  These curves are sometimes known as
\emph{Carnot-Carath\'eodory geodesics} and they have been extensively studied in general
sub-Riemannian manifolds, see \cite{montgomery}.  By the aforementioned regularity result any
length-minimizing curve in $(\stres,d)$ is a geodesic.  In this section we follow the approach
in \cite[\S~3]{stationary} to obtain a variational characterization of the geodesics.

Let $\ga:I\to\stres$ be a $C^2$ horizontal curve.  A \emph{smooth
variation} of $\ga$ is a $C^2$ map $F:I\times J\to\stres$, where $J$
is an open interval around the origin, such that $F(s,0)=\ga(s)$.  We
denote $\ga_\eps(s)=F(s,\eps)$.  Let $X_{\eps}(s)$ be the vector field
along $\ga_{\eps}$ given by $(\ptl F/\ptl\eps)(s,\eps)$.  Trivially
$[X_{\eps},\dot{\ga}_{\eps}]=0$.  Let $X=X_{0}$.  We say that the
variation is \emph{admissible} if the curves $\ga_\eps$ are horizontal
and have fixed extreme points.  For such a variation the vector field
$X$ vanishes at the endpoints of $\ga$ and satisfies
\[
0=\dot{\ga}\big(\escpr{X,V}\big)-2\,\escpr{X_{h},J(\dot{\ga})}.
\]
The equation above characterizes the vector fields along $\ga$
associated to admissible variations.  By using the first variation of
length in Riemannian geometry we can prove the following result, see
\cite[Proposition 3.1]{stationary} for details.

\begin{proposition}
\label{prop:eqgeo}
Let $\ga:I\to\stres$ be a $C^2$ horizontal curve parameterized by
arc-length.  Then $\ga$ is a critical point of length for any
admissible variation if and only if there is $\lambda\in\rr$ such that
$\ga$ satisfies the second order ordinary differential equation
\begin{equation}
\label{eq:geodesic}
D_{\dot{\ga}}\dot{\ga}+2\lambda\,J(\dot{\ga})=0.
\end{equation}
\end{proposition}

We will say that a $C^2$ horizontal curve $\ga$ is a \emph{geodesic of
curvature} $\la$ in $(\stres, g_h)$ if $\ga$ is parameterized by
arc-length and satisfies equation~\eqref{eq:geodesic}.  Observe that
the parameter $\la$ in \eqref{eq:geodesic} changes to $-\la$ for the
reversed curve $\ga(-s)$, while it is preserved for the antipodal
curve $-\ga(s)$.  In general, any isometry of $(\stres,g)$ preserving
the horizontal distribution transforms geodesics in geodesics since it
respects the connection $D$ of $g$ and commutes with $J$.

Given a point $p\in\stres$, a unit horizontal vector $v\in
T_{p}\stres$, and $\lambda\in\rr$, we denote by $\ga_{p,v}^\lambda$
the unique solution to \eqref{eq:geodesic} with initial conditions
$\ga(0)=p$ and $\dot{\ga}(0)=v$.  The curve $\ga_{p,v}^\lambda$ is a
geodesic since it is horizontal and parameterized by arc-length (the
functions $\escpr{\dot{\ga},V}$ and $|\dot{\ga}|^2$ are constant along
any solution of \eqref{eq:geodesic}).  Clearly for any right
translation $R_{q}$ we have $R_{q}(\ga^\la_{p,v})=\ga^\la_{p\cdot
q,v\cdot q}$.

Now we compute the geodesics in Euclidean coordinates.  Consider a
$C^2$ smooth curve $\ga=(x_{1},y_{1},x_{2},y_{2})\in\stres$
parameterized by arc-length $s$.  We denote
$\ddot{\ga}=(\ddot{x}_{1},\ddot{y}_{1},\ddot{x}_{2},\ddot{y}_{2})$.
The tangent and normal projections of $\ddot{\ga}$ onto $T\stres$ and
$(T\stres)^\bot$ are given respectively by $D_{\dot{\ga}}\dot{\ga}$
and $\text{II}(\dot{\ga},\dot{\ga})\,\eta$, where $\text{II}$ is the
second fundamental form of $\stres$ in $\rr^4$ with respect to the
unit normal vector $\eta(p)=p$.  Hence we obtain
\begin{equation}
\label{eq:ddotga}
\ddot{\ga}=D_{\dot{\ga}}\dot{\ga}-\ga.
\end{equation}
As a consequence equation \eqref{eq:geodesic} reads
\[
\ddot{\ga}+\ga+2\la\,(i\cdot\dot{\ga})=0.
\]
If we denote $z_{n}=x_{n} +iy_{n}$ then the previous equation is
equivalent to
\[
\ddot{z}_{n}+z_{n}+2\la i\,\dot{z}_{n}=0,\qquad n=1,2.
\]
Therefore, an explicit integration gives for $n=1,2$
\begin{equation}
\label{eq:integrate}
z_{n}(s)=C_{1n}\,\exp\{(-\la+\sqrt{1+\la^2})\,is\}+
C_{2n}\,\exp\{-(\la+\sqrt{1+\la^2})\,is\},
\end{equation}
where $C_{1n}$ and $C_{2n}$ are complex constants. Thus, if we
denote $C^r_{mn}=\text{Re}(C_{mn})$ and $C^i_{mn}=\text{Im}(C_{mn})$
then we have
\begin{align*}
x_{n}(s)&=(C^r_{1n}+C^r_{2n})\,\cos(\la
s)\,\cos(\sqrt{1+\la^2}\,s)+(C^r_{1n}-C^r_{2n})\,\sin(\la
s)\,\sin(\sqrt{1+\la^2}\,s)
\\
&+(C^i_{1n}+C^i_{2n})\,\sin(\la
s)\,\cos(\sqrt{1+\la^2}\,s)+(C^i_{2n}-C^i_{1n})\,\cos(\la
s)\,\sin(\sqrt{1+\la^2}\,s),
\\
y_{n}(s)&=(C^i_{1n}+C^i_{2n})\,\cos(\la
s)\,\cos(\sqrt{1+\la^2}\,s)-(C^i_{2n}-C^i_{1n})\,\sin(\la
s)\,\sin(\sqrt{1+\la^2}\,s)
\\
&-(C^r_{1n}+C^r_{2n})\,\sin(\la
s)\,\cos(\sqrt{1+\la^2}\,s)+(C^r_{1n}-C^r_{2n})\,\cos(\la
s)\,\sin(\sqrt{1+\la^2}\,s).
\end{align*}

Suppose that $\ga(0)=(x_{1}^0,y_{1}^0,x_{2}^0,y_{2}^0)$ and
$\dot{\ga}(0)=(u_{1}^0,w_{1}^0,u_{2}^0,w_{2}^0)$. It is easy to see
from \eqref{eq:integrate} that
\begin{align*}
C^r_{1n}+C^r_{2n}&=x^0_{n},\qquad
C^r_{1n}-C^r_{2n}=\frac{w^0_{n}+\la x^0_{n}}{\sqrt{1+\la^2}},
\\
C^i_{1n}+C^i_{2n}&=y^0_{n},\qquad \,
C^i_{2n}-C^i_{1n}=\frac{u^0_{n}-\la y^0_{n}}{\sqrt{1+\la^2}}.
\end{align*}
So, by substituting the previous equalities in the expressions of
$x_{n}(s)$ and $y_{n}(s)$ we obtain
\begin{align}
\label{eq:geocoor}
x_{n}(s)&=x^0_{n}\,\cos(\la
s)\,\cos(\sqrt{1+\la^2}\,s)+\frac{w^0_{n}+\la
x^0_{n}}{\sqrt{1+\la^2}}\,\sin(\la s)\,\sin(\sqrt{1+\la^2}\,s)
\\
\notag
&+y^0_{n}\,\sin(\la s)\,\cos(\sqrt{1+\la^2}\,s)+\frac{u^0_{n}-\la
y^0_{n}}{\sqrt{1+\la^2}}\,\cos(\la s)\,\sin(\sqrt{1+\la^2}\,s).
\\
\notag
y_{n}(s)&=y^0_{n}\,\cos(\la s)\,\cos(\sqrt{1+\la^2}\,s)-
\frac{u^0_{n}-\la y^0_{n}}{\sqrt{1+\la^2}}\,\sin(\la s)\,
\sin(\sqrt{1+\la^2}\,s)
\\
\notag
&-x^0_{n}\,\sin(\la s)\,\cos(\sqrt{1+\la^2}\,s)+\frac{w^0_{n}+\la
x^0_{n}}{\sqrt{1+\la^2}}\,\cos(\la s)\,\sin(\sqrt{1+\la^2}\,s).
\end{align}
We conclude that the geodesic $\ga^\la_{p,v}$ is given for any
$s\in\rr$ by
\begin{align}
\label{eq:geoconjunta}
\ga^\la_{p,v}(s)&=\cos(\la s)\,\cos(\sqrt{1+\la^2}\,s)\,p+
\frac{\sin(\la s)\,\sin(\sqrt{1+\la^2}\,s)}{\sqrt{1+\la^2}}\,(\la
p-J(v))
\\
\notag
&-\sin(\la s)\,\cos(\sqrt{1+\la^2}\,s)\,V(p)+
\frac{\cos(\la s)\,\sin(\sqrt{1+\la^2}\,s)}{\sqrt{1+\la^2}}\,(\la
V(p)+v).
\end{align}
In particular, for $\la=0$ we get
\[
\ga^0_{p,v}(s)=\cos(s)\,p+\sin(s)\,v,
\]
which is a horizontal great circle of $\stres$.  This was already
observed in \cite[Lemma~7.1]{chmy}.

Now we prove a characterization of the geodesics that will be useful in
Section~\ref{sec:classification}.  The result also shows that the geodesics are horizontal
lifts via the Hopf fibration $\mathcal{F}:\stres\to\sph^2$ of the geodesic circles in $\sph^2$,
see \cite[Thm.~1.26]{montgomery} for a general statement for principal bundles.

\begin{lemma}
\label{lem:geofunction}
Let $\ga:I\to\stres$ be a $C^2$ horizontal curve parameterized by
arc-length. The following assertions are equivalent
\begin{itemize}
\item[(i)] $\ga$ is a geodesic of curvature $\la$ in $(\stres,g_h)$,
\item[(ii)] $\escpr{\ddot{\ga},J(\dot{\ga})}=-2\la$,
\item[(iii)] the Hopf fibration $\mathcal{F}(\ga)$ is a piece of a
geodesic circle in $\sph^2$ with constant geodesic curvature $\la$
in $\sph^2$.
\end{itemize}
\end{lemma}

\begin{proof}
As $\ga$ is horizontal and parameterized by arc-length we have
\begin{align*}
0&=\dot{\ga}\,(\escpr{\dot{\ga},\dot{\ga}})=2\,\escpr{D_{\dot{\ga}}
\dot{\ga},
\dot{\ga}},
\\
0&=\dot{\ga}\,(\escpr{\dot{\ga},V(\ga)})=
\escpr{D_{\dot{\ga}}\dot{\ga},V(\ga)}+\escpr{\dot{\ga},J(\dot{\ga})}=
\escpr{D_{\dot{\ga}}\dot{\ga},V(\ga)}.
\end{align*}
As $\{\dot{\ga},J(\dot{\ga}),V(\ga)\}$ is an orthonormal basis of
$T\stres$ along $\ga$, we deduce that $D_{\dot{\ga}}\dot{\ga}$
is proportional to $J(\dot{\ga})$ at any point of $\ga$. On the other
hand from \eqref{eq:ddotga} we have
\[
\escpr{D_{\dot{\ga}}\dot{\ga},J(\dot{\ga})}=\escpr{\ddot{\ga}
+\ga,J(\dot{\ga})}=\escpr{\ddot{\ga},J(\dot{\ga})},
\]
where in the second equality we have used that the position vector
field $\eta(p)=p$ in $\rr^4$ provides a unit normal to $\stres$.  This
proves that (i) and (ii) are equivalent.

Let us see that (i) is equivalent to (iii).  Note that
$\mathcal{F}(R_{q}(p))=(L_{\overline{q}}\circ R_{q})(\mathcal{F}(p))$
for any $p,q\in\stres$.  Hence we only have to prove the claim for a
geodesic $\ga$ leaving from $p=(1,0,0,0)$.  Let
$v=(\cos\theta)\,E_{1}(p)+(\sin\theta)\,E_{2}(p)$ be the initial
velocity of such a geodesic.  A direct computation from
\eqref{eq:geocoor} shows that the Euclidean coordinates
$(y_1,x_2,y_2)$ of the curve $c=\mathcal{F}(\ga)$ are given~by
\begin{align*}
y_1(s)&=1-\frac{2}{1+\la^2}\,\sin^2(\sqrt{1+\la^2}\,s),
\\
\notag
x_2(s)&=\frac{-\sin(2\sqrt{1+\la^2}\,s)}{\sqrt{1+\la^2}}\,\sin\theta+
\frac{2\la\,\sin^2(\sqrt{1+\la^2}\,s)}{1+\la^2}\,\cos\theta.
\\
\notag y_2(s)&=\frac{\sin(2\sqrt{1+\la^2}\,s)}
{\sqrt{1+\la^2}}\,\cos\theta+\frac{2\la\,\sin^2(\sqrt{1+\la^2}\,s)}
{1+\la^2}\,\sin\theta.
\end{align*}
From the equations above it is not difficult to check that the
binormal vector to $c$ in $\rr^3$ is
$|\dot{c}\wedge\ddot{c}|^{-1}(\dot{c}\wedge\ddot{c})(s)=
(1+\la^2)^{-1/2}\,(\la,\sin\theta,\cos\theta)$.  It follows that the
curve $c$ lies inside a Euclidean plane and so, it must be a piece of
a geodesic circle in $\sph^2$.  Moreover, the geodesic curvature of
$c$ in $\sph^2$ with respect to the unit normal vector given by
$|c\wedge\dot{c}|^{-1}\,(c\wedge\dot{c})$ equals $\la$.  This proves
that (i) implies (iii).  Conversely, let us suppose that
$c=\mathcal{F}(\ga)$ is a piece of a geodesic circle of curvature
$\la$ in $\sph^2$.  We consider the geodesic $\ga^{\la}_{p,v}$ in
$(\stres,g_{h})$ with initial conditions $p=\ga(0)$ and
$v=\dot{\ga}(0)$.  The previous arguments and the uniqueness of
constant geodesic curvature curves in $\sph^2$ for given initial
conditions imply that $\mathcal{F}(\ga^\la_{p,v})=c$.  By using the
uniqueness of the horizontal lifts of a curve we conclude that
$\ga=\ga^\la_{p,v}$.
\end{proof}

In the next result we show that the topological behaviour of a
geodesic in $(\stres,g_{h})$ depends on the curvature of the geodesic.
Recall that $\mathcal{T}_{\rho}$ denotes the Clifford torus consisting
of the pairs $(z_{1},z_{2})\in\stres$ such that $|z_{1}|^2=\rho^2$.

\begin{proposition}
\label{prop:geoclosed}
Let $\ga:\rr\to\stres$ be a complete geodesic of curvature $\la$.
Then $\ga$ is a closed curve diffeomorphic to a circle if and only if
$\la/\sqrt{1+\la^2}$ is a rational number.  Otherwise $\ga$ is
diffeomorphic to a straight line and there is a right translation
$R_{q}$ such that $R_{q}(\ga)$ is a dense subset inside a Clifford
torus $\mathcal{T}_{\rho}$.
\end{proposition}

\begin{proof}
In order to characterize when $\ga$ is a closed curve diffeomorphic to
a circle it would be enough to analyze the equality
$\ga(s_{1})=\ga(s_{2})$ from \eqref{eq:geoconjunta}.  However we will
prove the proposition by using the description of a geodesic contained
inside a Clifford torus $\mathcal{T}_{\rho}$.

We shall use complex notation for the points in $\stres$.  Let
$q=(z_{1},z_{2})\in\mathcal{T}_{\rho}$.  It is easy to check that there are only two unit
horizontal vectors in $T_{q}\mathcal{T}_{\rho}$.  These are $w=i\cdot(\alpha z_{1},-\alpha^{-1}
z_{2})$ and $-w$, where $\alpha=\rho^{-1}\sqrt{1-\rho^2}$.  Take the geodesic
$\ga^\la_{q,w}=(z_{1}(s),z_{2}(s))$ of curvature $\la$.  A direct computation from
\eqref{eq:geocoor} gives us
\[
|z_{1}(s)|^2=\rho^2\left (\cos^2(\sqrt{1+\la^2}\,s)+
\frac{(\la+\alpha)^2}{1+\la^2}\,\sin^2(\sqrt{1+\la^2}\,s)\right ),
\]
so that $\ga^\la_{q,w}$ is entirely contained in $\mathcal{T}_{\rho}$
if and only if $\la=(2\rho^2-1)/(2\rho\sqrt{1-\rho^2})$.  Consider the
map $\varphi(x,y)=(\rho\,\exp(2\pi ix),\sqrt{1-\rho^2}\,\exp(2\pi
iy))$, which is a diffeomorphism between the flat torus
$\rr^2/\mathbb{Z}^2$ and $\mathcal{T}_{\rho}$ .  If we choose the
curvature $\la$ as above and we put $q=\varphi(\theta,\theta')$ then
we deduce from \eqref{eq:geocoor} that
\[
\ga^\la_{q,w}(s)=\varphi\left(\frac{(\sqrt{1+\la^2}-\la)\,s}{2\pi}+\theta,
\frac{-(\la+\sqrt{1+\la^2})\,s}{2\pi}+\theta'\right).
\]
This implies that $\ga^\la_{q,w}$ is a reparameterization of
$\varphi(r(t))$, where $r(t)=mt+n$ is a straight line in
$\rr^2/\mathbb{Z}^2$ with slope
\[
m=\frac{\la+\sqrt{1+\la^2}}{\la-\sqrt{1+\la^2}}=
\frac{(\la/\sqrt{1+\la^2})+1}{(\la/\sqrt{1+\la^2})-1}.
\]
As a consequence $\ga^\la_{q,w}$ is a closed curve diffeomorphic to a
circle if and only if $\la/\sqrt{1+\la^2}$ is a rational number.
Otherwise $\ga^\la_{q,w}$ is a dense curve in $\mathcal{T}_{\rho}$
diffeomorphic to a straight line.

Finally, let us consider any complete geodesic $\ga=\ga^\la_{p,v}$ in
$(\stres,g_h)$.  After applying a right translation we can suppose
that $p=(1,0)$ and $v=(0,\exp(i\theta))$.  Let $\rho\in (0,1)$ so that
$\la/\sqrt{1+\la^2}=2\rho^2-1$.  Take the point
$q=(\rho,\sqrt{1-\rho^2}\,i\exp(i\theta))\in\mathcal{T}_{\rho}$.  It
is easy to check that the vector $v\cdot q$ coincides with the unit
horizontal vector $w\in T_{q}\mathcal{T}_{\rho}$ such that
$\ga^\la_{q,w}\sub\mathcal{T}_{\rho}$.  The proof of the proposition
then follows by using that $R_{q}(\ga^\la_{p,v})=\ga^\la_{q,w}$ and
the properties previously shown for geodesics inside
$\mathcal{T}_\rho$.
\end{proof}

\begin{figure}[h]
\centering{\includegraphics[width=6cm]{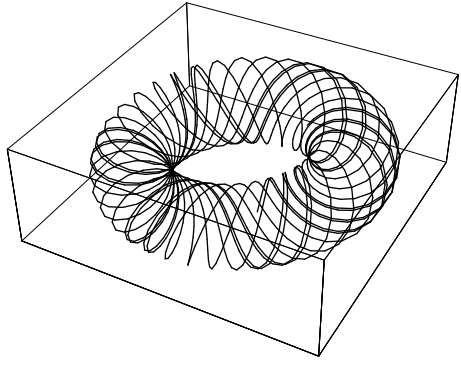}}
\caption{Stereographic projection from $\stres-\{\text{north pole}\}$
to $\rr^3$ of a sub-Riemannian geodesic which is dense inside a
Clifford torus.}
\label{fig:densegeo}
\end{figure}

We finish this section with some analytical properties for the vector
field associated to a variation of a curve which is a geodesic.  The
proofs use the same arguments as in Lemma 3.5 and Lemma 3.6 in
\cite{stationary}.

\begin{lemma}
\label{lem:jacobi}
Let $\ga:I\to\stres$ be a geodesic of curvature $\la$.  Suppose that
$X$ is the $C^1$ vector field associated to a variation of $\ga$ by
horizontal curves $\ga_\eps$ parameterized by arc-length.  Then we
have
\begin{itemize}
\item[(i)] The function $\la\,\escpr{X,V(\ga
)}+\escpr{X,\dot{\ga}}$ is
constant along $\ga$.
\item[(ii)] If any $\ga_\eps$ is a geodesic of curvature $\la$
and $X$ is $C^2$ smooth, then $X$ satisfies the second order differential
equation $D_{\dot{\ga}}D_{\dot{\ga}}{X}+R(X,\dot{\ga})\dot{\ga}+
2\la\,(J(D_{\dot{\ga}}{X})-\escpr{X,\dot{\ga}}\,V(\ga))=0$, where $R$
denotes the Riemannian curvature tensor in $(\stres,g)$.
\end{itemize}
\end{lemma}

The linear differential equation in Lemma~\ref{lem:jacobi} (ii) is the
\emph{Jacobi equation} for geodesics of curvature $\la$ in
$(\stres,g_{h})$.  We will call any solution of this equation a
\emph{Jacobi field} along $\ga$.

\section{Area-stationary surfaces with or without a volume constraint}
\label{sec:mainresult}
\setcounter{equation}{0}

In this section we introduce and characterize critical surfaces for
the area functional \eqref{eq:area} with or without a volume
constraint.  We also state without proof some properties for such
surfaces that will be useful in order to obtain classifications
results.  For a detailed development we refer the reader to \cite[\S
4]{stationary} and the references therein.

Let $\Sg\sub\stres$ be an oriented immersed surface of class $C^2$.
Consider a $C^1$ vector field $X$ with compact support on $\Sg$ and
tangent to $\stres$.  For $t$ small we denote
$\Sg_{t}=\{\exp_{p}(tX_{p});p\in\Sg\}$, which is an immersed surface.
Here $\exp_{p}$ is the exponential map of $(\stres,g)$ at the point
$p$.  The family $\{\Sg_t\}$, for $t$ small, is the \emph{variation}
of $\Sg$ induced by $X$.  Note that we allow the variations to move
the singular set $\Sg_{0}$ of $\Sg$.  Define $A(t):=A(\Sg_t)$.  If
$\Sg$ is the boundary of a region $\Om\sub\stres$ then we can consider
a $C^1$ family of regions $\Om_{t}$ such that $\Om_{0}=\Om$ and
$\ptl\Om_t=\Sg_t$.  We define $V(t):=V(\Om_t)$.  We say that the
variation induced by $X$ is \textit{volume-preserving} if $V(t)$ is
constant for any $t$ small enough.  We say that $\Sg$ is
\emph{area-stationary} if $A'(0)=0$ for any variation of $\Sg$.  In
case that $\Sg$ encloses a region $\Om$, we say that $\Sg$ is
\emph{area-stationary under a volume constraint} or
\emph{volume-preserving area-stationary} if $A'(0)=0$ for any
volume-preserving variation of $\Sg$.

Suppose that $\Om$ is the region bounded by a $C^2$ embedded compact
surface $\Sg$.  We shall always choose the unit normal $N$ to $\Sg$ in
$(\stres,g)$ pointing into $\Om$.  The computation of $V'(0)$ is well
known, and it is given by (\cite[\S 9]{simon})
\begin{equation}
\label{eq:1stvol}
V'(0)=\int_{\Om}\divv X\,dv=-\int_\Sg u\,d\Sg,
\end{equation}
where $u=\escpr{X,N}$.  It follows that $u$ has mean zero whenever the
variation is volume-preserving.  Conversely, it was proved in
\cite[Lemma 2.2]{bdce} that, given a $C^1$ function $u:\Sg\to\rr$ with
mean zero, we can construct a volume-preserving variation of $\Om$ so
that the normal component of $X$ equals $u$.

\begin{remark}
\label{re:bdce} For a compact immersed $C^2$ surface $\Sg$ in $\stres$
there is a notion of volume enclosed by $\Sg$.  The first variation
for this volume functional is given by \eqref{eq:1stvol}.  We refer
the reader to \cite[p.~125]{bdce} for details.
\end{remark}

Now assume that the divergence relative to $\Sg$ of the horizontal
Gauss map $\nu_h$ defined in \eqref{eq:nuh} satisfies
$\divv_\Sg\nu_h\in L^1(\Sg)$.  In this case the first variation of the
area functional $A(t)$ can be obtained as in \cite[Lemma
4.3]{stationary}.  We get
\begin{equation}
\label{eq:1stvar}
A'(0)=\int_\Sg
u\,\big(\divv_\Sg\nu_h\big)\,d\Sg
-\int_\Sg\divv_\Sg\big(u\,(\nu_h)^\top\big)\,d\Sg,
\end{equation}
where $(\nu_{h})^\top$ is the projection of $\nu_h$ onto the tangent
space to $\Sg$.

Let $\Sg$ be a $C^2$ immersed surface in $\stres$ with a $C^1$ unit
normal vector $N$.  Outside the singular set $\Sg_{0}$ of $\Sg$ we
define the \emph{mean curvature} $H$ in $(\stres,g_h)$ by the equality
\begin{equation}
\label{eq:mc}
-2H(p):=(\divv_{\Sg}\nu_{h})(p),\qquad p\in\Sg-\Sg_{0}.
\end{equation}
This notion of mean curvature agrees with the ones introduced in
\cite{chmy} and \cite{hp}.  We say that $\Sg$ is a \emph{minimal
surface} if $H\equiv 0$ on $\Sg-\Sg_{0}$.  By using variations
supported in $\Sg-\Sg_{0}$, the first variation of area
\eqref{eq:1stvar}, and the first variation of volume
\eqref{eq:1stvol}, we deduce that the mean curvature of $\Sg-\Sg_{0}$
is respectively zero or constant if $\Sg$ is area-stationary or
volume-preserving area-stationary.  In $\Sg-\Sg_{0}$ we can consider
the orthonormal basis $\{Z,S\}$ defined in \eqref{eq:nuh} and
\eqref{eq:ese}, so that we get from \eqref{eq:mc}
\[
-2H=\escpr{D_{Z}\nu_{h},Z}+\escpr{D_{S}\nu_{h},S}.
\]
It is easy to check (\cite[Lemma 4.2]{stationary}) that for any tangent
vector $X$ to $\Sg$ we have
\[
D_{X}\nu_{h}=|N_{h}|^{-1}\,\big(\escpr{D_{X}N,Z}+\escpr{N,V}\,
\escpr{X,\nu_{h}}\big)\,Z+\escpr{Z,X}\,V.
\]
In particular by taking $X=Z$ and $X=S$ we deduce the following
expression for the mean curvature
\begin{equation}
\label{eq:mc2}
2H=|N_{h}|^{-1}\,\text{II}(Z,Z),
\end{equation}
where $\text{II}$ is the second fundamental form of $\Sg$ with respect
to $N$ in $(\stres,g)$.

On the other hand, by the arguments in \cite[Thm.~4.8]{stationary},
any characteristic curve $\ga$ of a $C^2$ immersed surface $\Sg$
satisfies
\begin{equation}
\label{eq:geoh} D_{\dot{\ga}}\dot{\ga}=-2H\,J(\dot{\ga}).
\end{equation}
From the previous equality we deduce that $\Sg-\Sg_{0}$ is a ruled
surface in $(\stres,g_{h})$ whenever $H$ is constant, see also
\cite[Cor.~6.10]{hp}.

\begin{theorem}
\label{th:ruled}
Let $\Sg$ be an oriented $C^2$ immersed surface in $(\stres,g_{h})$
with constant mean curvature $H$ outside the singular set.  Then any
characteristic curve of $\Sg$ is an open arc of a geodesic of
curvature $H$ in $(\stres,g_{h})$.
\end{theorem}

Now we describe the configuration of the singular set $\Sg_{0}$ of a
constant mean curvature surface $\Sg$ in $(\stres,g_{h})$.  The set
$\Sg_0$ was studied by J.-H.~Cheng, J.-F.~Hwang, A.~Malchiodi and
P.~Yang \cite{chmy} for surfaces with bounded mean curvature inside
the first Heisenberg group.  As indicated by the authors in
\cite[Lemma 7.3]{chmy} and \cite[Proof of Thm.~E]{chmy}, their local
arguments also apply for spherical pseudohermitian $3$-manifolds.  We
gather their results in the following theorem.

\begin{theorem}[{\cite[Theorem~B]{chmy}}]
\label{th:chmy}
Let $\Sg$ be a $C^2$ oriented immersed surface in $(\stres,g_h)$ with
constant mean curvature $H$ off of the singular set $\Sg_{0}$.  Then
$\Sg_{0}$ consists of isolated points and $C^1$ curves with
non-vanishing tangent vector.  Moreover, we have
\begin{itemize}
\item[(i)] $($\cite[Thm.~3.10]{chmy}$)$ If $p\in\Sg_{0}$ is
isolated then there exists $r>0$ and $\la\in\rr$ with $|\la|=|H|$ such
that the set described as
\[
D_{r}(p)=\{\gamma_{p,v}^\la(s);v\in T_p\Sg,\,|v|=1,\,s\in [0,r)\},
\]
is an open neighborhood of $p$ in $\Sg$.  \vspace{0,1cm} \item[(ii)]
$($\cite[Prop.~3.5 and Cor.~3.6]{chmy}$)$ If $p$ is contained in a
$C^1$ curve $\Ga\subset\Sg_{0}$ then there is a neighborhood $B$ of
$p$ in $\Sg$ such that $B\cap\Gamma$ is a connected curve and $B-\Ga$
is the union of two disjoint connected open sets $B^+$ and $B^-$
contained in $\Sg-\Sg_0$.  Furthermore, for any $q\in\Ga\cap B$ there
are exactly two geodesics $\ga_{1} ^\la\sub B^+$ and $\ga_{2}^\la\sub
B^-$ starting from $q$ and meeting transversally $\Ga$ at $q$ with
opposite initial velocities.  The curvature $\la$ does not depend on
$q\in\Ga\cap B$ and satisfies $|\la|=|H|$.
\end{itemize}
\end{theorem}

\begin{remark}
\label{rem:lambdaorientation}
The relation between $\lambda$ and $H$ depends on the value of the
normal $N$ to $\Sg$ in the singular point $p$.  If $N_{p}=V_{p}$ then
$\lambda=H$, whereas $\lambda=-H$ when $N_{p}=-V_{p}$.  In case
$\la=H$ the geodesics $\ga^\la$ in Theorem~\ref{th:chmy} are
characteristic curves of $\Sg$.
\end{remark}

The characterization of area-stationary surfaces with or without a
volume constraint in $(\stres,g_{h})$ is similar to the one obtained by
M.~Ritor\'e and the second author in \cite[Thm.~4.16]{stationary}.
We can also improve, as in \cite[Prop.  4.19]{stationary}, the $C^1$
regularity of the singular curves of an area-stationary surface.

\begin{theorem}
\label{th:constant}
Let $\Sg$ be an oriented $C^2$ immersed surface in $\stres$.  The
followings assertions are equivalent
\begin{itemize}
\item[(i)] $\Sg$ is area-stationary $($resp.  volume-preserving
area-stationary$)$ in $(\stres,g_h)$.
\item[(ii)] The mean curvature of $\Sg-\Sg_{0}$ is
zero $($resp.  constant$)$ and the characteristic curves meet
orthogonally the singular curves when they exist.
\end{itemize}
Moreover, if $(\emph{i})$ holds then the singular curves of $\Sg$
are $C^2$ smooth.
\end{theorem}

\begin{example}
\label{ex:syt}
1.  Every totally geodesic $2$-sphere in $(\stres,g)$ is a compact
minimal surface in $(\stres, g_h)$.  In fact, for any $q\in\stres$,
the $2$-sphere $\stres\cap q^\bot$ is the union of all the points
$\ga^0_{p,v}(s)$ where $p=-i\cdot q$, the unit vector $v\in T_p\stres$
is horizontal, and $s\in [0,\pi]$.  These spheres have two singular
points at $p$ and $-p$.  In particular they are area-stationary
surfaces by Theorem~\ref{th:constant}.

2.  For any $\rho\in (0,1)$ the Clifford torus $\mathcal{T}_{\rho} $ has no singular points
since the vertical vector $V$ is tangent to this surface.  We consider the unit normal vector
to $\mathcal{T}_{\rho}$ in $(\stres,g)$ given for $q=(z_{1},z_{2})$ by $N(q)=(\alpha
z_{1},-\alpha^{-1} z_{2})$, where $\alpha=\rho^{-1}\,\sqrt{1-\rho^2}$.  As $\escpr{N,V}=0$ then
we have $N=N_{h}=\nu_{h}$ and so $Z=J(N)$.  Let $\la=(2\rho^2-1)/(2\rho\sqrt{1-\rho^2})$.  It
was shown in the proof of Proposition~\ref{prop:geoclosed} that the geodesic $\ga^\la_{q,w}$
with $w=Z(q)$ is entirely contained in $\mathcal{T}_{\rho}$.  The tangent vector to this
geodesic equals $Z$ since the singular set is empty.  We conclude that $\ga^\la_{q,w}$ is a
characteristic curve of $\mathcal{T}_{\rho}$.  By using \eqref{eq:geoh} we deduce that
$\mathcal{T}_{\rho}$ has constant mean curvature $H=(2\rho^2-1)/(2\rho\sqrt{1-\rho^2})$ with
respect to the normal $N$. By Theorem~\ref{th:constant} the surface $\mathcal{T}_{\rho}$ is
volume-preserving area-stationary for any $\rho\in (0,1)$. Moreover, $\mathcal{T}_\rho$ is
area-stationary for $\rho=\sqrt{2}/2$.
\end{example}

The previous examples were found in \cite{chmy}.  In
\cite[Theorem~E]{chmy}, J.-~H.~Cheng, J.-F. Hwang, A.~Malchiodi and
P.~Yang described the possible topological types for a compact surface
with bounded mean curvature inside a spherical pseudohermitian
$3$-manifold were.  More precisely, they proved the following result.

\begin{theorem}
\label{th:genus}
Let $\Sg$ be an immersed $C^2$ compact, connected, oriented surface in
$(\stres,g_h)$ with bounded mean curvature outside the singular set.
If $\Sg$ contains an isolated singular point then $\Sg$ is
homeomorphic to a sphere.  Otherwise $\Sg$ is homeomorphic to a torus.
\end{theorem}

\section{Classification results for complete stationary surfaces}
\label{sec:classification}

An immersed surface $\Sg\subset\stres$ is \emph{complete} if it is
complete in $(\stres,g)$.  We say that a complete, noncompact,
oriented $C^2$ surface $\Sg$ is volume-preserving area-stationary if
it has constant mean curvature off of the singular set and the
characteristic curves meet orthogonally the singular curves when they
exist.  By Theorem \ref{th:constant} this implies that $\Sg$ is a
critical point for the area functional of any variation with compact
support of $\Sg$ such that the ``volume enclosed'' by the perturbed
region is constant, see Remark~\ref{re:bdce}.

\subsection{Complete surfaces with isolated singularities}
\label{subsec:isolated}

It was shown in \cite[Proof of Cor.~F]{chmy} that any $C^2$ compact,
connected, embedded, minimal surface in $(\stres,g_{h})$ with an
isolated singular point coincides with a totally geodesic $2$-sphere
in $(\stres,g)$.  In this section we generalize this result for
complete immersed surfaces with constant mean curvature.  First we
describe the surface which results when we join two certain points in
$\stres$ by all the geodesics of the same curvature.

For $p=(1,0,0,0)$ and $\la\in\rr$, let $\ga_{\theta}$ be the geodesic
of curvature $\la$ in $(\stres,g_{h})$ with initial conditions
$\ga_{\theta}(0)=p$ and
$\dot{\ga}_{\theta}(0)=v=(\cos\theta)\,E_{1}(p)+(\sin\theta)\,E_{2}(p)$.
By \eqref{eq:geocoor} the Euclidean coordinates of $\ga_{\theta}$ are
given by
\begin{align}
\label{eq:geocoore1}
x_{1}(s)&=\cos(\la s)\,\cos(\sqrt{1+\la^2}\,s)+
\frac{\la}{\sqrt{1+\la^2}}\,\sin(\la s)\,\sin(\sqrt{1+\la^2}\,s),
\\
\notag
y_{1}(s)&=-\sin(\la s)\,\cos(\sqrt{1+\la^2}\,s)+
\frac{\la}{\sqrt{1+\la^2}}\,\cos(\la s)\,\sin(\sqrt{1+\la^2}\,s),
\\
\notag
x_{2}(\theta,s)&=\frac{1}{\sqrt{1+\la^2}}\,\sin(\sqrt{1+\la^2}\,s)\,
\cos(\theta-\la s),
\\
\notag
y_{2}(\theta,s)&=\frac{1}{\sqrt{1+\la^2}}\,\sin(\sqrt{1+\la^2}\,s)\,
\sin(\theta-\la s).
\end{align}
We remark that the functions $x_{1}(s)$ and $y_{1}(s)$ in
\eqref{eq:geocoore1} do not depend on $\theta$.  We define
$\mathcal{S}_{\la}$ to be the set of points $\ga_{\theta}(s)$ where
$\theta\in [0,2\pi]$ and $s\in [0,\pi/\sqrt{1+\la^2}]$.  From
\eqref{eq:geocoore1} it is clear that the point
$p_{\la}:=\ga_{\theta}(\pi/\sqrt{1+\la^2})$ is the same for any
$\theta$.  In fact, we have
\[
p_{\la}=-\cos\left(\frac{\la\pi}{\sqrt{1+\la^2}}\right)\,p+
\sin\left(\frac{\la\pi}{\sqrt{1+\la^2}}\right)\,V(p).
\]
It follows that $p_{\la}$ moves along the vertical great circle of
$\stres$ passing through $p$.  Note that $p_{0}=-p$ and $p_{\la}\to p$
when $\la\to\pm\infty$.  We will call $p$ and $p_{\la}$ the
\emph{poles} of~ $\mathcal{S}_{\la}$.  Observe that $\mathcal{S}_{0}$
coincides with a totally geodesic $2$-sphere in $(\stres,g)$, see
Example~\ref{ex:syt}.  From \eqref{eq:geocoore1} we also see that
$\mathcal{S}_{\la}$ is invariant under any rotation $r_{\theta}$ in
\eqref{eq:rtheta}.

\begin{proposition}
\label{prop:spheres}
The set $\mathcal{S}_{\la}$ is a $C^2$ embedded volume-preserving
area-stationary $2$-sphere with constant mean curvature
$\la$ off of the poles.
\end{proposition}

\begin{proof}
We consider the $C^\infty$ map $F:[0,2\pi]\times
[0,\pi/\sqrt{1+\la^2}]\to\stres$ defined by
$F(\theta,s)=\ga_{\theta}(s)$.  Clearly $F(0,s)=F(2\pi,s)$,
$F(\theta,0)=p$ and $F(\theta,\pi/\sqrt{1+\la^2})=p_{\la}$.  Suppose
that $F(\theta_{1} ,s_{1})= F(\theta_{2},s_{2})$ for $\theta_{i}\in
[0,2\pi)$ and $s_{i}\in (0,\pi/\sqrt{1+\la^2})$.  This is equivalent
to that $\ga_{\theta_{1}}(s_{1})=\ga_{\theta_{2}}(s_{2})$.  For
$\la\neq 0$ the function $y_{1}(s)$ in \eqref{eq:geocoore1} is
monotonic on $(0,\pi/\sqrt{1+\la^2})$ since its first derivative
equals $(1+\la^2)^{-1/2}\,\sin(\la s)\,\sin(\sqrt{1+\la^2}\,s)$.  For
$\la=0$ we have $x_1(s)=\cos(s)$, which is decreasing on $(0,\pi)$.
So, equality $\ga_{\theta_{1}}(s_{1})=\ga_{\theta_{2}}(s_{2})$ implies
$s_{1}=s_{2}=s_{0}$.  Moreover, the equalities between the
$x_{2}$-coordinates and the $y_{2}$-coordinates of
$\ga_{\theta_{1}}(s_{0})$ and $\ga_{\theta_{2}}(s_{0})$ yield
$\theta_{1} =\theta_{2}$.  The previous arguments show that
$\mathcal{S}_{\la}$ is homeomorphic to a $2$-sphere.

Note that $(\ptl F/\ptl s)(\theta,s)=\dot{\ga}_{\theta}(s)$, which is
a horizontal vector.  Let $X_{\theta}(s):=(\ptl
F/\ptl\theta)(\theta,s)$.  By Lemma~\ref{lem:jacobi} (ii) this is a
Jacobi vector field along $\ga_\theta$ vanishing for $s=0$ and
$s=\pi/\sqrt{1+\la^2}$.  The components of $X_{\theta}$ with respect
to $\dot{\ga}_{\theta}$ and $V(\ga_{\theta})$ can be computed from
\eqref{eq:geocoore1} so that we get
\begin{align*}
\escpr{X_{\theta}(s),\dot{\ga}_{\theta}(s)}&=\left(\frac{\ptl
x_{2}}{\ptl\theta}\,\frac{\ptl x_{2}}{\ptl s}+\frac{\ptl
y_{2}}{\ptl\theta}\,\frac{\ptl y_{2}}{\ptl
s}\right)(\theta,s)=-\frac{\la\,\sin^2(\sqrt{1+\la^2}\,s)}{1+\la^2},
\\
\escpr{X_{\theta}(s),V(\ga_{\theta}(s))}&=\left(x_{2}\,\frac{\ptl
y_{2}}{\ptl\theta}-y_{2}\,\frac{\ptl
x_{2}}{\ptl\theta}\right) (\theta,s)=\frac{\sin^2(\sqrt{1+\la^2}\,s)}{1+\la^2}.
\end{align*}
It follows that $X_{\theta}(s)$ has a non-trivial vertical component
for $s\in (0,\pi/\sqrt{1+\la^2})$.  As a consequence,
$\mathcal{S}_{\la}$ with the poles removed is a $C^\infty$ smooth
embedded surface in $\stres$ without singular points.

To prove that $\mathcal{S}_{\la}$ is volume-preserving area-stationary
it suffices by Theorem~\ref{th:constant} to show that the mean
curvature is constant off of the poles.  Consider the unit normal
vector along $\mathcal{S}_{\la}-\{p,p_{\la}\}$ defined by
$N=(1-\escpr{X_{\theta},\dot{\ga}_{\theta}}^2)^{-1/2}\,
(-\escpr{X_{\theta},V(\ga_\theta)}\,J(\dot{\ga}_{\theta})+\escpr{X_{\theta}
,J(\dot{\ga}_{\theta})}\,V(\ga_\theta))$.  The characteristic vector
field associated to $N$ is given by
$Z(\theta,s)=\dot{\ga}_{\theta}(s)$.  By using \eqref{eq:geoh} we
deduce that $\mathcal{S}_{\la}-\{p,p_{\la}\}$ has constant mean
curvature $\la$ with respect to $N$.  To complete the proof it is
enough to observe that $\mathcal{S}_{\la}$ is also a $C^2$ embedded
surface around the poles.  This is a consequence of
Remark~\ref{re:graphs} below.
\end{proof}

\begin{remark}
\label{re:graphs}
The surface $\mathcal{S}_{\la}$ can be described as the union of two
radial graphs over the $x_{2}y_{2}$ plane.  Let
$r_{\la}=1/\sqrt{1+\la^2}$ and $\var(r)=r_{\la}\,\arcsin(r/r_{\la})$
for $r\in [0,r_{\la}]$.  We can see from \eqref{eq:geocoore1} that the
lower half of $\mathcal{S}_{\la}$ is given by
\begin{align*}
x_{1}(r)&=\sqrt{1-(r/r_{\la})^2}\,\cos(\la\var(r))+\la r\,\sin(\la\var(r)),
\\
y_{1}(r)&=\la r\,\cos(\la\var(r))-\sqrt{1-(r/r_{\la})^2}\,\sin(\la\var(r)),
\end{align*}
where $r=(x_{2}^2+y_{2}^2)^{1/2}$ belongs to $[0,r_{\la}]$.
Similarly, the upper half of $\mathcal{S}_{\la}$ can be described as
\begin{align*}
x_{1}(r)&=-\sqrt{1-(r/r_{\la})^2}\,\cos(\la\psi(r))+\la r\,\sin(\la\psi(r)),
\\
y_{1}(r)&=\la r\,\cos(\la\psi(r))+\sqrt{1-(r/r_{\la})^2}\,\sin(\la\psi(r)),
\end{align*}
where $\psi(r)=\pi r_{\la}-\var(r)$.  The poles are the points
obtained for $r=0$ and they are singular points of
$\mathcal{S}_{\la}$.  From the equations above it can be shown that
$\mathcal{S}_{\la}$ is $C^2$ around these points.  Moreover,
$\mathcal{S}_{\la}$ is $C^3$ around the north pole if and only if
$\la=0$, i.e., $\mathcal{S}_{\la}$ is a totally geodesic $2$-sphere in
$(\stres,g)$.
\end{remark}

Now we can prove our first classification result.

\begin{theorem}
\label{th:spheres}
Let $\Sg$ be a complete, connected, oriented, immersed $C^2$ surface
with constant mean curvature in $(\stres,g_{h})$ .  If $\Sg$ contains
an isolated singular point then $\Sg$ is congruent with a
sphere~$\mathcal{S}_{\la}$.
\end{theorem}

\begin{proof}
We reproduce the arguments in \cite[Thm.~6.1]{stationary}.  Let $H$ be
the mean curvature of $\Sg$ with respect to a unit normal vector $N$.
After a right translation of $\stres$ we can assume that $\Sg$ has an
isolated singularity at $p=(1,0,0,0)$.  Suppose that $N_p=V(p)$.  By
Theorem~\ref{th:chmy} (i) and Remark~\ref{rem:lambdaorientation},
there exists a neighborhood $D_{r}$ of $p$ in $\Sg$ which consists of
all the geodesics of curvature $\la=H$ and length $r$ leaving from
$p$.  By using Theorem~\ref{th:ruled} and the completeness of $\Sg$ we
deduce that these geodesics can be extended until they meet a singular
point.  As $\Sg$ is immersed and connected we conclude that
$\Sg=\mathcal{S}_{\la}$.  Finally, if $N_p=-V(p)$ we repeat the
previous arguments by using geodesics of curvature $\la=-H$ and we
obtain that $\Sg=\phi(\mathcal{S}_{\la})$, where $\phi$ is the
isometry of $(\stres,g)$ given by
$\phi(x_{1},y_{1},x_{2},y_{2})=(x_{1},-y_{1},x_{2},-y_{2})$.  Clearly
$\phi$ preserves the horizontal distribution so that $\Sg$ is
congruent with $\mathcal{S}_\la$.
\end{proof}

\subsection{Complete surfaces with singular curves}
\label{subsec:singularcurves}

In this section we follow the arguments in \cite[\S 6]{stationary} to
describe complete area-stationary surfaces in $\stres$ with or without
a volume constraint and non-empty singular set consisting of $C^2$
curves.  For such a surface we know by Theorem~\ref{th:constant} that
the characteristic curves meet orthogonally the singular curves.
Moreover, if the surface is compact then it is homeomorphic to a torus
by virtue of Theorem~\ref{th:genus}.

We first study in more detail the behaviour of the characteristic
curves of a volume-preserving area-stationary surface far away from a
singular curve.  Let $\Ga:I\to\stres$ be a $C^2$ curve defined on an
open interval.  We suppose that $\Ga$ is horizontal with arc-length
parameter $\eps\in I$.  We denote by $\ddot{\Ga}$ the covariant
derivative of $\dot{\Ga}$ for the flat connection on $\rr^4$.  Note
that $\{\Ga,\dot{\Ga},J(\dot{\Ga}),V(\Ga)\}$ is an orthonormal basis
of $\rr^4$ for any $\eps\in I$.  Thus we get
\begin{equation}
\label{eq:gamma2}
\ddot{\Ga}=-\Ga+h\,J(\dot{\Ga}),
\end{equation}
where $h=\escpr{\ddot{\Ga},J(\dot{\Ga})}$.  Fix $\la\in\rr$.  For any
$\eps\in I$, let $\ga_{\eps}(s)$ be the geodesic in $(\stres,g_h)$ of
curvature $\la$ with initial conditions $\ga_{\eps}(0)=\Ga(\eps)$ and
$\dot{\ga}_{\eps}(0)=J(\dot{\Ga}(\eps))$.  Clearly $\ga_{\eps}$ is
orthogonal to $\Ga$ at $s=0$.  By equation \eqref{eq:geoconjunta} we
have
\begin{align}
\label{eq:geocoor2}
\ga_{\eps}(s)&=\left(\cos(\la
s)\,\cos(\sqrt{1+\la^2}\,s)+\frac{\la\,\sin(\la
s)\,\sin(\sqrt{1+\la^2}\,s)}{\sqrt{1+\la^2}}\right)\,\Ga(\eps)
\\
\notag
&+\frac{\sin(\la s)\,\sin(\sqrt{1+\la^2}\,s)}{\sqrt{1+\la^2}}\,\,
\dot{\Ga}(\eps)+\frac{\cos(\la s)\,
\sin(\sqrt{1+\la^2}\,s)}{\sqrt{1+\la^2}}\,\,J(\dot{\Ga}(\eps))
\\
\notag
&+\left(-\sin(\la s)\,\cos(\sqrt{1+\la^2}\,s)+
\frac{\la\,\cos(\la s)\,\sin(\sqrt{1+\la^2}\,s)}{\sqrt{1+\la^2}}\right)
\,V(\Ga(\eps)).
\end{align}
We define the $C^1$ map $F(\eps,s)=\ga_{\eps}(s)$, for $\eps\in I$ and
$s\in[0,\pi/\sqrt{1+\la^2}]$.  Note that $(\ptl F/\ptl
s)(\eps,s)=\dot{\ga}_{\eps}(s)$.  We define $X_{\eps}(s):=(\ptl
F/\ptl\eps)(\eps,s)$.  In the next result we prove some properties
of~$X_{\eps}$.

\begin{lemma}
\label{lem:jacobi3}
In the situation above, $X_{\eps}$ is a Jacobi vector field along
$\ga_{\eps}$ with $X_{\eps}(0)=\dot{\Ga}(\eps)$.  For any $\eps\in I$
there is a unique $s_{\eps}\in (0,\pi/\sqrt{1+\la^2})$ such that
$\escpr{X_{\eps}(s_{\eps}),V(\ga_{\eps}(s_{\eps}))}=0$.  We have
$\escpr{X_{\eps},V(\ga_{\eps})}<0$ on $(0,s_{\eps})$ and
$\escpr{X_{\eps},V(\ga_{\eps})}>0$ on $(s_{\eps},\pi/\sqrt{1+\la^2})$.
Moreover $X_{\eps}(s_{\eps})=J(\dot{\ga}_{\eps}(s_{\eps}))$.
\end{lemma}

\begin{proof}
We denote by $a(s)$, $b(s)$, $c(s)$ and $d(s)$ the components of
$\ga_{\eps}(s)$ with respect to the orthonormal basis
$\{\Ga,\dot{\Ga},J(\dot{\Ga}),V(\Ga)\}$, see \eqref{eq:geocoor2}.
By using \eqref{eq:gamma2} we have that
\begin{align*}
\frac{d}{d\eps}\,J(\dot{\Ga}(\eps))&=\frac{d}{d\eps}\,(i\cdot\dot{\Ga}(\eps))
=i\cdot\ddot{\Ga}(\eps)=-V(\Ga(\eps))-h(\eps)\dot{\Ga}(\eps),
\\
\frac{d}{d\eps}\,V(\Ga(\eps))&=\frac{d}{d\eps}\,(i\cdot\Ga(\eps))
=i\cdot\dot{\Ga}(\eps)=J(\dot{\Ga}(\eps)).
\end{align*}
From here, the definition of $X_{\eps}$, and \eqref{eq:gamma2} we obtain
\[
X_{\eps}(s)=-b(s)\,\Ga(\eps)+(a(s)-h(\eps)c(s))\,\dot{\Ga}(\eps)
+(d(s)+h(\eps)b(s))\,J(\dot{\Ga}(\eps))-c(s)\,V(\Ga(\eps)).
\]
It follows that $X_{\eps}(0)=\dot{\Ga}(\eps)$ and that $X_{\eps}$ is a
$C^\infty$ vector field along $\ga_{\eps}$.  Moreover, $X_{\eps}$ is a
Jacobi vector field along $\ga_{\eps}$ by Lemma~\ref{lem:jacobi} (ii).
The vertical component of $X_{\eps}$ can be computed from
\eqref{eq:geocoor2} so that we get
\begin{align*}
\escpr{X_{\eps},V(\ga_{\eps})}(s)&=\escpr{X_{\eps}(s),i\cdot\ga_{\eps}(s)}=
2\,(b(s)d(s)-a(s)c(s))+h(\eps)\,(b(s)^2+c(s)^2)
\\
&=\frac{\sin(\sqrt{1+\la^2}\,s)}{\sqrt{1+\la^2}}\,
\left(\frac{\sin(\sqrt{1+\la^2}\,s)}{\sqrt{1+\la^2}}\,h(\eps)
-2\cos(\sqrt{1+\la^2}\,s)\right).
\end{align*}
Thus $\escpr{X_{\eps}(s_{\eps}),V(\ga_{\eps}(s_{\eps}))}=0$ for some
$s_{\eps}\in (0,\pi/\sqrt{1+\la^2})$ if and only if
\begin{equation}
\label{eq:despeje}
h(\eps)=2\sqrt{1+\la^2}\,\cot(\sqrt{1+\la^2}\,s_{\eps} ).
\end{equation}
From \eqref{eq:despeje} we obtain the existence and uniqueness of
$s_{\eps}$ as that as the sign of $\escpr{X_{\eps},V(\ga_\eps)}$.

Now we use Lemma \ref{lem:jacobi} (i) and that
$X_{\eps}(0)=\dot{\Ga}(\eps)$ to deduce that the function given by
$\la\,\escpr{X_{\eps},V(\ga_\eps)}+\escpr{X_{\eps},\dot{\ga}_{\eps}}$
vanishes along $\ga_{\eps}$.  In particular, $X_{\eps}(s_{\eps})$ is a
horizontal vector orthogonal to $\dot{\ga}_{\eps}(s_{\eps})$.
Finally, a straightforward computation gives us
\begin{align*}
\escpr{X_{\eps},J(\dot{\ga}_{\eps})}(s)&=b(s)\dot{d}(s)
-(a(s)-h(\eps)c(s))\,\dot{c}(s)+(d(s)+h(\eps)b(s))\,\dot{b}(s)-\dot{a}(s)c(s)
\\
\notag
&=\frac{\sin(2\sqrt{1+\la^2}\,s)}{2\sqrt{1+\la^2}}\,h(\eps)
-\cos(2\sqrt{1+\la^2}\,s),\quad s\in[0,\pi/\sqrt{1+\la^2}].
\end{align*}
By using \eqref{eq:despeje} we see that the expression above equals
$1$ for $s=s_{\eps}$.  This completes the proof.
\end{proof}

In the next result we construct immersed surfaces with constant mean
curvature boun\-ded by two singular curves.  Geometrically we only
have to leave from a given horizontal curve by segments of orthogonal
geodesics of the same curvature.  The length of these segments is
indicated by the \emph{cut function} $s_{\eps}$ defined in
Lemma~\ref{lem:jacobi3}.  We also characterize when the resulting
surfaces are area-stationary with or without a volume constraint.

\begin{proposition}
\label{prop:sigmala}
Let $\Ga$ be a $C^{k+1}$ $(k\geq 1)$ horizontal curve in $\stres$
parameterized by arc-length $\eps\in I$.  Consider the map $F:I\times
[0,\pi/\sqrt{1+\la^2}]\to\stres$ defined by $F(\eps,s)=\ga_{\eps}(s)$,
where $\ga_{\eps}$ is the geodesic of curvature $\la$ with initial
conditions $\Ga(\eps)$ and $J(\dot{\Ga}(\eps))$.  Let $s_{\eps}$ be
the function introduced in Lemma~\ref{lem:jacobi3}, and let
$\Sg_{\la}(\Ga):=\{F(\eps,s);\,\eps\in I, \,s\in [0,s_{\eps}]\}$.  Then
we have
\begin{itemize}
\item[(i)] $\Sg_{\la}(\Ga)$ is an immersed surface of class $C^k$ in
$\stres$.
\item[(ii)] The singular set of $\Sg_\la(\Ga)$ consists of
two curves $\Ga(\eps)$ and $\Ga_{1}(\eps):=F(\eps,s_{\eps})$.
\item[(iii)] There is a $C^{k-1}$ unit normal vector $N$ to
$\Sg_{\la}(\Ga)$ in $(\stres,g)$ such that $N=V$ on $\Ga$ and
$N~=-V$ on $\Ga_{1}$.
\item[(iv)] The curve $\ga_{\eps}(s)$ for
$s\in(0,s_{\eps})$ is a characteristic curve of $\Sg_{\la}(\Ga)$ for
any $\eps\in I$.  In particular, if $k\geq 2$ then $\Sg_{\la}(\Ga)$ has
constant mean curvature $\la$ in $(\stres,g_h)$ with respect to~$N$.
\item[(v)] If\/ $\Ga_{1}$ is a $C^2$ smooth curve then the geodesics
$\ga_{\eps}$ meet orthogonally $\Ga_{1}$ if and only if $s_{\eps}$ is
constant along $\Ga$.  This condition is equivalent to that $\Ga$ is a
geodesic in $(\stres,g_h)$.
\end{itemize}
\end{proposition}

\begin{proof}
That $F$ is a $C^k$ map is a consequence of \eqref{eq:geocoor2} and
the fact that $\Ga$ is $C^{k+1}$.  Consider the vector fields $(\ptl
F/\ptl\eps)(\eps,s)=X_{\eps}(s)$ and $(\ptl F/\ptl
s)(\eps,s)=\dot{\ga}_{\eps}(s)$.  By Lemma \ref{lem:jacobi3} we deduce
that the differential of $F$ has rank two for any $(s,\eps)\in I\times
[0,\pi/\sqrt{1+\la^2})$, and that the tangent plane to
$\Sg_{\la}(\Ga)$ is horizontal only for the points in $\Ga$ and
$\Ga_{1}$.  This proves (i) and (ii).

Consider the $C^{k-1}$ unit normal vector to the immersion $F:I\times
[0,\pi/\sqrt{1+\la^2})\to~\stres$ given by
$N=(1-\escpr{X_{\eps},\dot{\ga}_{\eps}}^2)^{-1/2}\,
(\escpr{X_{\eps},V(\ga_\eps)}\,J(\dot{\ga}_{\eps})-
\escpr{X_{\eps},J(\dot{\ga}_{\eps})}\,V(\ga_\eps))$.  Since we have
$X_{\eps}(0)=\dot{\Ga}(\eps)$ and
$X_{\eps}(s_{\eps})=J(\dot{\ga}_{\eps}(s_{\eps}))$ it follows that
$N=V$ along $\Ga$ and $N=-V$ along $\Ga_{1}$.  On the other hand, the
characteristic vector field associated to $N$ is
\[
Z(\eps,s)=-\frac{\escpr{X_{\eps}(s),V(\ga_\eps(s))}}
{|\escpr{X_{\eps}(s),V(\ga_\eps(s))}|}\,\,\dot{\ga}_{\eps}(s),
\qquad\eps\in I, \ \ s\neq 0, s_{\eps},
\]
and so $Z(\eps,s)= \dot{\ga_{\eps}}(s)$ whenever $s\in (0,s_{\eps})$ by
Lemma~\ref{lem:jacobi3}.  This fact and \eqref{eq:geoh} prove (iv).

Finally, suppose that $\Ga_{1}$ is a $C^2$ smooth curve.  In this
case, the cut function $s(\eps)=s_{\eps}$ is $C^1$, and the tangent
vector to $\Ga_{1}$ is given by
\[
\dot{\Ga}_{1}(\eps)=X_{\eps}(s_{\eps})+\dot{s}(\eps)\,
\dot{\ga}_{\eps}(s_{\eps}).
\]
As $X_\eps(s_{\eps})=J(\dot{\ga}_{\eps}(s_{\eps}))$ we conclude that
the geodesics $\ga_\eps$ meet $\Ga_1$ orthogonally if and only if
$s(\eps)$ is a constant function.  By \eqref{eq:despeje} the function
$h=\escpr{\ddot{\Ga},J(\dot{\Ga})}$ is constant along $\Ga$.  By
Lemma~\ref{lem:geofunction} this is equivalent to that $\Ga$ is a
geodesic.
\end{proof}

\begin{remark}
\label{re:reverse}
1.  In the proof of Proposition~\ref{prop:sigmala} we have shown that
if we extend the surface $\Sg_\la(\Ga)$ by the geodesics $\ga_\eps$
beyond the singular curve $\Ga_1$ then the resulting surface has mean
curvature $-\la$ beyond $\Ga_{1}$.  As indicated in
Theorem~\ref{th:chmy} (ii), to obtain an extension of $\Sg_{\la}(\Ga)$
with constant mean curvature $\la$ we must leave from $\Ga_{1}$ by
geodesics of curvature $-\lambda$.

2.  Let $\Ga:I\to\stres$ be a $C^{k+1}$ $(k\geq 1)$ horizontal curve
parameterized by arc-length.  We consider the geodesic
$\widetilde{\ga}_{\eps}$ of curvature $\la$ and initial conditions
$\Ga(\eps)$ and $-J(\dot{\Ga}(\eps))$.  By following the arguments in
Lemma \ref{lem:jacobi3} and Proposition \ref{prop:sigmala} we can
construct the surface
$\widetilde{\Sg}_{\la}(\Ga):=\{\widetilde{\ga}_{\eps}(s);\ \eps\in
I,\, s\in [0,\widetilde{s}_{\eps}]\}$, which is bounded by two
singular curves $\Ga$ and $\Ga_2$.  The value $\widetilde{s}_{\eps}$
is defined as the unique $s\in (0,\pi/\sqrt{1+\la^2})$ such that
$\escpr{\widetilde{X}_{\eps}, V(\widetilde{\ga}_\eps)}(s)=0$.  Here
$\widetilde{X}_{\eps}$ is the Jacobi vector field associated to the
variation $\{\widetilde{\ga}_{\eps}\}$.  The cut function
$\widetilde{s}_{\eps}$ satisfies the equality
\begin{equation}
\label{eq:despeje2}
h(\eps)=-2\sqrt{1+\la^2}\,\cot(\sqrt{1+\la^2}\,\,\widetilde{s}_{\eps}),
\end{equation}
where $h=\escpr{\ddot{\Ga},J(\dot{\Ga})}$.  From \eqref{eq:despeje} it
follows that $s_{\eps}+\widetilde{s}_{\eps} =\pi/\sqrt{1+\la^2}$.  The
vector $\widetilde{X}_{\eps}$ coincides with
$-J(\dot{\widetilde{\ga}}_{\eps})$ for $s=\widetilde{s}_{\eps}$.  We
can define a unit normal $\widetilde{N}$ satisfying $\widetilde{N}=V$
on $\Ga$ and $\widetilde{N}=-V$ on $\Ga_2$.  For $k\ge 2$ we deduce
that $\Sg_{\la}(\Ga)\cup\widetilde{\Sg}_{\la}(\Ga)$ is an oriented
immersed surface with constant mean curvature $\la$ outside the
singular set and at most three singular curves.
\end{remark}

Now we shall use Proposition \ref{prop:sigmala} and
Remark~\ref{re:reverse} to obtain examples of complete surfaces with
constant mean curvature outside a non-empty set of singular curves.
Taking into account Theorem~\ref{th:constant} and
Proposition~\ref{prop:sigmala} (v), if we also require the surfaces to
be volume-preserving area-stationary then the initial curve $\Ga$ must
be a geodesic.

\begin{example}[The torus $\mathcal{C}_{0,\la}$]
\label{ex:singulartori}
Let $\Ga$ be the horizontal great circle of $\stres$ parameterized by
$\Ga(\eps)=(\cos(\eps),0,\sin(\eps),0)$ (the geodesic of curvature
$\mu=0$ with initial conditions $p=(1,0,0,0)$ and $v=E_1(p)$).  For
any $\la\in\rr$ let $\mathcal{C}_{0,\la}$ be the union of the surfaces
$\Sg_{\la}(\Ga)$ and $\widetilde{\Sg}_{\la}(\Ga)$ introduced in
Proposition \ref{prop:sigmala} and Remark \ref{re:reverse}.  The
resulting surface is $C^\infty$ outside the singular set and has
constant mean curvature $\la$.  The cut functions $s_\eps$ and
$\widetilde{s}_\eps$ associated to $\Ga$ can be obtained from
\eqref{eq:despeje} and $\eqref{eq:despeje2}$, so that we get
$s_\eps=\widetilde{s}_\eps=\pi/(2\sqrt{1+\la^2})$.  By using
\eqref{eq:geocoor2} we can compute the map $F(\eps,s)=\ga_\eps(s)$
defined for $\eps\in [0,2\pi]$ and $s\in [0,\pi/(2\sqrt{1+\la^2})]$.
In particular we can give an explicit expression for the singular
curve $\Ga_1(\eps)$, which is a horizontal great circle different from
$\Ga$.  Let $\eps_0\in (0,\pi)$ such that $\cot(\eps_0)=-\la$.  It is
easy to check that $\Ga_1(\eps_0)=\exp(i\theta_1)\cdot p$ and
$\dot{\Ga}_1(\eps_0)=\exp(i\theta_1)\cdot v$, where
$\theta_1=3\pi/2-(\la\pi)/(2\sqrt{1+\la^2})$.  By using the uniqueness
of the geodesics we deduce that
$\Ga_1(\eps+\eps_0)=\exp(i\theta_1)\cdot\Ga(\eps)$.  With similar
arguments we obtain that
$\Ga_2(\eps+\widetilde{\eps}_0)=\exp(i\theta_2)\cdot\Ga(\eps)$, where
$\widetilde{\eps}_0=\pi-\eps_0$ and $\theta_2=\theta_1-\pi$.  Note
that $\exp(i\theta_1)\cdot p=-\exp(i\theta_2)\cdot p$.  As any great
circle of $\stres$ is invariant under the antipodal map $q\mapsto -q$,
we conclude that $\Ga_1$ and $\Ga_2$ are different parameterizations
of the same horizontal great circle.

\begin{figure}[h]
\centering{\includegraphics[height=6cm]{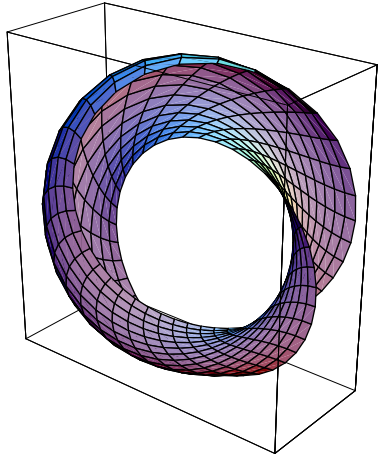}}
\caption{Stereographic projection from $\stres-\{\text{north pole}\}$
to $\rr^3$ of one half of the surface $\mathcal{C}_{0,\la}$.  It
consists of the union of all the geodesics of curvature $\la$ and
length $\pi/(2\sqrt{1+\la^2})$ connecting two singular circles.}
\label{fig:halfsingulartorus}
\end{figure}

In Figure~\ref{fig:halfsingulartorus} we see that the surface $\Sg_\la(\Ga)$ is embedded.  To
prove this note that the function $(x_1y_1+x_2y_2)(\eps,s)$ only depends on $s$, and its first
derivative with respect to $s$ equals $(1+\la^2)^{-1/2}\,\sin(2\la
s)\,\sin(2\sqrt{1+\la^2}\,s)$, which does not change sign on $(0,\pi/(2\sqrt{1+\la^2}))$.  Thus
if $F(\eps_1,s_1)=F(\eps_2,s_2)$ for some $\eps_i\in [0,2\pi)$ and $s_i\in
[0,\pi/(2\sqrt{1+\la^2})]$ then $s_1=s_2$, which clearly implies $\eps_1=\eps_2$.  Similarly we
obtain that $\widetilde{\Sg}_\la(\Ga)$ is embedded.  On the other hand, observe that $2(x_1
y_2-x_2y_1)(\eps,s)=\sin(2\sqrt{1+\la^2}\,s)/\sqrt{1+\la^2}$ on $\Sg_\la(\Ga)$, whereas the
same function evaluated on $\widetilde{\Sg}_\la$ equals
$-\sin(2\sqrt{1+\la^2}\,s)/\sqrt{1+\la^2}$.  It follows that $\mathcal{C}_{0,\la}$ is an
embedded surface outside the singular curves.  Finally, a long but easy computation shows that
there is a system of coordinates $(u_1,u_2,u_3,u_4)$ such that $\mathcal{C}_{0,\la}$ can be
expressed as union of certain graphs $u_1=f_i(u_2,u_3)$ and $u_4=g_i(u_2,u_3)$, $i=1,2$,
defined over an annulus of the $u_2u_3$-plane.  The functions $f_i$ and $g_i$ are $C^2$ near
the singular curves.  This proves that $\mathcal{C}_{0,\la}$ is a volume-preserving
area-stationary embedded torus with two singular curves.
\end{example}

\begin{figure}[h]
\centering{\includegraphics[width=5cm]{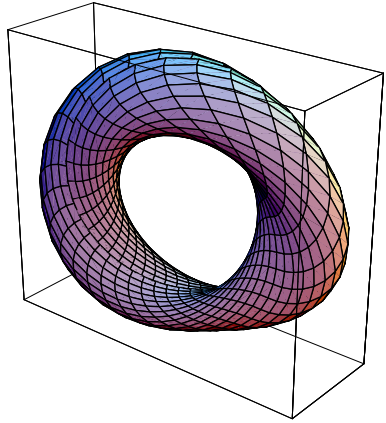}}
\caption{Stereographic projection from $\stres-\{\text{north pole}\}$
to $\rr^3$ of an embedded torus $\mathcal{C}_{0,\la}$.}
\label{fig:singulartorus}
\end{figure}

\begin{example}[The surfaces $\mathcal{C}_{\mu,\la}$]
\label{ex:smula}
Let $\Ga$ be the geodesic of curvature $\mu$ in $(\stres,g_h)$ with
initial conditions $p=(1,0,0,0)$ and $v=E_1(p)$.  We know that the
function $h=\escpr{\ddot{\Ga},J(\dot{\Ga})}$ equals $-2\mu$ along
$\Ga$ by Lemma~\ref{lem:geofunction}.  For any $\la\in\rr$ we consider
the union $\Sg_\la(\Ga)\cup\widetilde{\Sg}_\la(\Ga)$, which is a
$C^\infty$ surface with constant mean curvature $\la$ outside the
singular curves $\Ga$, $\Ga_1$ and $\Ga_2$.  By using
Lemma~\ref{lem:geofunction} (ii) we can prove that any $\Ga_i$ is a
geodesic of curvature $\mu$.  The cut functions $s_\eps$ and
$\widetilde{s}_{\eps}$ are determined by equalities \eqref{eq:despeje}
and \eqref{eq:despeje2}.  Define $\eps_\mu$ as the unique
$\eps\in(0,\pi/\sqrt{1+\mu^2})$ such that
$\cot(\sqrt{1+\mu^2}\,\eps_\mu)=-\la/\sqrt{1+\mu^2}$.  Let
$\widetilde{\eps}_\mu=\pi/\sqrt{1+\mu^2}-\eps_\mu$.  Easy computations
from \eqref{eq:geocoor2} show that
\begin{eqnarray*}
&&\Gamma_1(s_\mu)=\exp(i\theta_1)\cdot p,\quad
\dot{\Gamma}_1(s_\mu)=\exp(i\theta_1)\cdot v,\\
&&\Gamma_2(\tilde{s}_\mu)=\exp(i\theta_2)\cdot p,\quad
\dot{\Gamma}_2(\tilde{s}_\mu)=\exp(i\theta_2)\cdot v,
\end{eqnarray*}
where $\theta_1=3\pi/2-\lambda s_\eps-\mu \eps_\mu$ and
$\theta_2=\pi/2-\la\widetilde{s}_\eps-\mu\widetilde{\eps}_\mu$.  By
the uniqueness of the geodesics we deduce that
$\Gamma_1(\eps+\eps_\mu)=\exp(i\theta_1)\cdot \Gamma(\eps)$ and
$\Gamma_2(\eps+\widetilde{\eps}_\mu)=\exp(i\theta_2)\cdot
\Gamma(\eps)$.  In general $\Ga_1\neq \Ga_2$ so that we can extend the
surface by geodesics orthogonal to $\Ga_i$ of the same curvature.  As
we pointed out in Remark~\ref{re:reverse} and according with the
initial velocity of $\Ga_i$, in order to preserve the constant mean
curvature $\la$ we must consider the surfaces
$\widetilde{\Sg}_{-\la}(\Ga_1)$ and $\Sg_{-\la}(\Ga_2)$.  Two new
singular curves $\Ga_{12}$ and $\Ga_{22}$ are obtained.  It is
straightforward to check that, after a translation of the parameter
$\eps$, we have $\Ga_{12}=\exp(i\theta_{12})\cdot \Gamma$ and
$\Ga_{22}=\exp(i\theta_{22})\cdot \Gamma$, where
$\theta_{12}=\theta_1+\pi/2+\lambda \tilde{s}_\eps-\mu \eps_\mu$ and
$\theta_{22}=\theta_2+ 3\pi/2 +\lambda s_\eps-\mu \tilde{\eps}_\mu$.
Let $\tilde{\theta}_1=\theta_{12}-\theta_1$ and
$\tilde{\theta}_2=\theta_{22}-\theta_2$.  We repeat this process by
induction so that at any step $k+1$ we leave from the singular curves
$\Ga_{1k}$ and $\Ga_{2k}$ by the corresponding orthogonal geodesics of
curvature $(-1)^k\la$.  We denote by $\mathcal{C}_{\mu,\la}$ the union
of all these surfaces.  After a translation of $\eps$, any singular
curve $\Ga_{jk}$ is of the form $\exp(i\theta_{jk})\cdot\Ga$, where
the angles are given by $\theta_{j\, 2m}=
m(\theta_j+\tilde{\theta}_j)$ and $\theta_{j\, 2m+1}= (m+1)\theta_j+
m\tilde{\theta}_j$.  This implies that all the singular curves are
geodesics of curvature $\mu$ and their projections to $\sph^2$ via the
Hopf fibration give the same geodesic circle.  It follows by
uniqueness of the horizontal lifts that two singular curves meeting at
one point must coincide as subsets of $\stres$.  In fact it is
possible that two singular curves coincide.  For example, the surface
$\mathcal{C}_{\mu,0}$ is a compact surface with two or four singular
curves (depending on if $\mu/\sqrt{1+\mu^2}$ is rational or not).  On
the other hand it can be shown that if $\mu/\sqrt{1+\mu^2}$ and
$\lambda/\sqrt{1+\lambda^2}$ are rational numbers and $(\lambda s_\eps
+\mu \eps_\mu)/\pi$ is irrational (take $\la=\mu=1/\sqrt{3}$) then
$\mathcal{C}_{\mu,\la}$ is a noncompact surface with infinitely many
singular curves.

The surface $\mathcal{C}_{\mu,\la}$ is $C^\infty$ off of the singular
set and has constant mean curvature $\la$.  A necessary condition to
get a surface which is also $C^2$ near the singular curves is that
$\Ga$ locally separates $\mathcal{C}_{\mu,\la}$ into two disjoint
domains, see Theorem~\ref{th:chmy} (ii).  By
Proposition~\ref{prop:geoclosed} this is equivalent to that
$\mu/\sqrt{1+\mu^2}$ is a rational number.  In such a case
$\mathcal{C}_{\mu,\la}$ is a volume-preserving area-stationary surface
by construction.  In general the surfaces $\mathcal{C}_{\mu,\la}$ are
not embedded.
\end{example}

Now we can classify complete area-stationary surfaces under a volume
constraint with a non-empty set of singular curves.

\begin{theorem}
\label{th:classification}
Let $\Sg$ be a complete, oriented, connected, $C^2$ immersed surface.
Suppose that $\Sg$ is volume-preserving area-stationary in $(\stres,g_{h})$ and $\Ga$ is a
connected singular curve of $\Sg$. Then $\Ga$ is a closed geodesic, and $\Sg$ is congruent with
a surface $\mathcal{C}_{\mu,\la}$.
\end{theorem}

\begin{proof}
By Theorem~\ref{th:constant} we have that $\Ga$ is a $C^2$ horizontal curve.  We can assume
that $\Ga$ is parameterized by arc-length.  We take the unit normal $N$ to $\Sg$ such that
$N=V$ along $\Ga$.  Let $H$ be the mean curvature of $\Sg$ with respect to $N$.  Let $p\in\Ga$.
By Theorem~\ref{th:chmy} (ii) and Remark~\ref{rem:lambdaorientation} there is a small
neighborhood $B$ of $p$ in $\Sg$ such that $B\cap\Ga$ is a connected curve separating $B$ into
two disjoint connected open sets foliated by geodesics $\ga_\eps$ of curvature $\la=H$ leaving
from $\Ga$.  These geodesics are characteristic curves of $\Sg$.  Moreover, by
Theorem~\ref{th:constant} they must leave from $\Ga$ orthogonally.  As $\Sg$ is complete and
connected we deduce that any $\ga_\eps$ can be extended until it meets a singular point.  Thus
there exists a small piece $\Ga'\sub\Ga$ containing $p$ and such that $\Sg_{\la}(\Ga')\sub\Sg$.
In particular we find another singular curve $\Ga_1'$ of $\Sg$ which is also $C^2$ smooth by
Theorem~\ref{th:constant}.  As $\Sg$ is volume-preserving area-stationary, any $\ga_\eps$ meet
$\Ga_1'$ orthogonally and so $\Ga'$ is a geodesic by Proposition~\ref{prop:sigmala} (v).  Since
$p\in\Ga$ is arbitrary we have proved that $\Ga$ is a geodesic in $(\stres,g_h)$.  That $\Ga$
is closed follows from Proposition~\ref{prop:geoclosed}; otherwise, the intersection of $\Ga$
with any open neighborhood of $p$ in $\Sg$ would have infinitely many connected components, a
contradiction with Theorem~\ref{th:chmy}~(ii). After applying a right translation $R_q$ and a
rotation $r_\theta$ we can suppose that $\Ga$ leaves from $p=(1,0,0,0)$ with velocity
$v=E_1(p)$.  By using again the local description of $\Sg$ around $\Ga$ in
Theorem~\ref{th:chmy} (ii) together with the completeness and the connectedness of $\Sg$, we
conclude that $\Sg$ is congruent with $\mathcal{C}_{\mu,\la}$.
\end{proof}

\subsection{Complete surfaces with empty singular set}
\label{subsec:empty}

Here we prove some classification results for complete constant mean
curvature surfaces with empty singular set.  Such a surface must be
area-stationary with or without a volume constraint by
Theorem~\ref{th:constant}.  Moreover, if the surface is compact then
it must be homeomorphic to a torus by Theorem~\ref{th:genus}.

The following result uses the behaviour of geodesics in $(\stres,g_h)$
described in Proposition~\ref{prop:geoclosed} to establish a strong
restriction on a compact embedded surface with constant mean
curvature.

\begin{theorem}
\label{th:irrat}
Let $\Sg$ be a $C^2$ compact, connected, embedded surface in
$(\stres,g_h)$ without singular points.  If $\Sg$ has constant mean
curvature $H$ such that $H/\sqrt{1+H^2}$ is an irrational number, then
$\Sg$ is congruent with a Clifford torus.
\end{theorem}

\begin{proof}
As $\Sg$ is compact with empty singular set we deduce by
Theorem~\ref{th:ruled} that there is a complete geodesic $\ga$ of
curvature $H$ contained in $\Sg$.  After a right translation $R_{q}$
we have, by Proposition~\ref{prop:geoclosed}, that $R_{q}(\ga)$ is a
dense subset of a Clifford torus $\mathcal{T}_{\rho}$.  By using that
$\Sg$ is compact, connected and embedded we conclude that
$R_{q}(\Sg)=\mathcal{T}_{\rho}$, proving the claim.
\end{proof}

In Remark~\ref{re:examplesrat} we will give examples showing that all
the hypotheses Theorem~\ref{th:irrat} are necessary.  We finish this
section with a characterization of the Clifford tori
$\mathcal{T}_\rho$ as the unique vertical surfaces with constant mean
curvature.  We say that a $C^1$ surface $\Sg\sub\stres$ is
\emph{vertical} if the vector field $V$ is tangent to $\Sg$.

\begin{proposition}
\label{prop:vertical}
Let $\Sg$ be a $C^2$ complete, connected, oriented, constant mean
curvature surface in $(\stres,g_{h})$.  If $\Sg$ is vertical then
$\Sg$ is congruent with a Clifford torus.
\end{proposition}

\begin{proof}
It is clear that $\Sg$ has no singular points.  Thus we can find by
Theorem~\ref{th:ruled} a complete geodesic $\ga$ contained in $\Sg$.
By Proposition~\ref{prop:geoclosed} there is a point $q\in\stres$ such
that $R_{q}(\ga)$ is contained inside a Clifford torus
$\mathcal{T}_{\rho}$.  By assumption, the vertical great circle
passing through any point of $\ga$ is entirely contained in
$R_{q}(\Sg)$.  Clearly the union of all these circles is
$\mathcal{T}_{\rho}$.  Finally as $\Sg$ is complete and connected we
conclude that $R_{q}(\Sg)=\mathcal{T}_{\rho}$.
\end{proof}

\section{Rotationally invariant constant mean curvature surfaces}
\label{sec:revolution}

In this section we classify $C^2$ constant mean curvature surfaces of
revolution in $(\stres,g_h)$.  We will follow arguments similar to
those in \cite[\S 5]{revolucion}.

Let $R$ be the great circle given by the intersection of $\stres$ with
the $x_{1}y_{1}$-plane.  The rotation $r_{\theta}$ of the
$x_{2}y_{2}$-plane defined in \eqref{eq:rtheta} is an isometry of
$(\stres,g)$ leaving invariant the horizontal distribution and fixing
$R$.  Let $\Sg$ be a $C^2$ surface in $\stres$ which is invariant
under any rotation $r_\theta$.  We denote by $\ga$ the generating
curve of $\Sg$ inside the hemisphere $\sph^2_+:=\{x_{2}\geq
0,y_{2}=0\}$.  If we parameterize $\ga=(x_{1},y_{1},x_{2})$ by
arc-length $s\in I$, then $\Sg-R$ is given in cylindrical coordinates
by $\phi(s,\theta)=r_{\theta}(\ga(s))=(x_{1}(s),y_{1}(s),x_{2}(s)
\cos\theta,x_{2}(s)\sin\theta)$.  Denote by $\{e_{1},e_{2}\}$ the
usual orthonormal frame in the Euclidean plane.  The tangent plane to
$\Sg-R$ is generated by the vector fields $\ptl_{1}:=e_1(\phi)$ and
$\ptl_{2}:=e_2(\phi)$.  Note that $|\ptl_{1}|=1$, $|\ptl_{2}|=x_2$ and
$\escpr{\ptl_{1},\ptl_{2}}=0$.  A unit normal vector along $\phi$ is
given by
\begin{equation}
\label{eq:normal}
N=(x_{2}\dot{y}_{1}-\dot{x}_{2} y_1,
x_{1}\dot{x}_{2}-\dot{x}_{1}x_{2},
(\dot{x}_{1}y_{1}-x_{1}\dot{y}_{1})\cos\theta,
(\dot{x}_{1}y_{1}-x_{1}\dot{y}_{1})\sin\theta).
\end{equation}
It follows that $|N_{h}|^2=\escpr{N,E_{1}}^2+\escpr{N,E_{2}}^2=
(\dot{x}_{1}y_{1}-x_{1}\dot{y}_1)^2+x_2^2$.  In particular, the
singular points of $\Sg$ are contained inside $R$.

Now we compute the mean curvature $H$ of $\Sg-R$ with respect to the
normal $N$ defined in \eqref{eq:normal}.  By equality \eqref{eq:mc2}
we know that $2H=|N_{h}|^{-1}\text{II}(Z,Z)$ and so, it is enough to
compute the second fundamental form $\text{II}$ of $\phi$ with respect
to $N$.  It is clear that the coefficients of $\text{II}$ in the basis
$\{\ptl_{1},\ptl_{2}\}$ are given by
$\text{II}_{ij}=\text{II}(\ptl_i,\ptl_j)=-\escpr{D_{\ptl_i}N,\ptl_j}=
\escpr{N,D_{e_i}\ptl_j}$.  On the other hand, if
$(a_{1j},a_{2j},a_{3j})$ are the coordinates of $\ptl_j$ in the
orthonormal basis $\{E_{1},E_{2},V\}$, then a straightforward
calculation by using \eqref{eq:christoffel} shows that the coordinates
of $D_{e_i}\ptl_j$ with respect to $\{E_{1},E_{2},V\}$ are
\[
\left(\frac{\ptl a_{1j}}{\ptl e_{i}}+a_{3i}\,a_{2j}-a_{2i}\,a_{3j},
\frac{\ptl a_{2j}}{\ptl e_{i}}-a_{3i}\,a_{1j}+a_{1i}\,a_{3j},
\frac{\ptl a_{3j}}{\ptl e_{i}}+a_{2i}\,a_{1j}-a_{1i}\,a_{2j}\right).
\]
This allows us to compute $\text{II}_{ij}$ and we obtain the following
\begin{align*}
\text{II}_{11}&=x_{2}\,(\ddot{x}_{1}\dot{y}_{1}-\dot{x}_{1}\ddot{y}_{1})
-\dot{x}_{2}\,(\ddot{x}_{1}y_{1}-x_{1}\ddot{y}_{1})
+\ddot{x}_{2}\,(\dot{x}_{1}y_{1}-x_{1}\dot{y}_{1}),
\\
\text{II}_{12}&=\text{II}_{21}=0,
\\
\text{II}_{22}&=x_{2}\,(x_{1}\dot{y}_{1}-\dot{x}_{1}y_{1}).
\end{align*}
On the other hand, the coordinates of the characteristic vector field
$Z$ with respect to $\{\ptl_{1},\ptl_{2}\}$ are $\escpr{Z,\ptl_{1}}$
and $x^{-2}_{2}\escpr{Z,\ptl_{2}}$.  Thus we can use equation
\eqref{eq:mc2} to deduce that the mean curvature of $\Sg-R$ with
respect to $N$ is
\[
2H=\frac{(x_{1}\dot{y}_{1}-\dot{x}_{1}y_{1})^3
+x_{2}^3\,\{x_{2}\,(\ddot{x}_{1}\dot{y}_{1}-\dot{x}_{1}\ddot{y}_{1})
-\dot{x}_{2}\,(\ddot{x}_{1}y_{1}-x_{1}\ddot{y}_{1})
+\ddot{x}_{2}\,(\dot{x}_{1}y_{1}-x_{1}\dot{y}_{1})\}} {x_{2}\,
(\dot{x}_{1}+\dot{y}_{1})^{3/2}}.
\]

Now we take spherical coordinates $(\omega,\tau)$ in
$\sph^2\equiv\{y_{2}=0\}$.  In precise terms, we choose $\omega\in
(-\pi/2,\pi/2)$ and $\tau\in\rr$ so that the Euclidean coordinates of
a point in $\sph^2$ different from the poles can be expressed as
$x_{1}=\cos\omega\,\cos\tau$, $y_{1}=\cos\omega\,\sin\tau$ and
$x_{2}=\sin\omega$.  The vector fields $\ptl_{\omega}$ and
$\ptl_{\tau}/(\cos\omega)$ provide an orthonormal basis of the tangent
plane to $\sph^2$ off of the poles.  The integral curves of
$\ptl_\omega$ and $\ptl_\tau$ are the meridians and the circles of
revolution about the $x_{2}$-axis, respectively.

Let $(\omega(s),\tau(s))$ with $\omega(s)\in [0,\pi/2)$ be the
spherical coordinates of the generating curve $\ga(s)$.  Denote by
$\sg(s)$ the oriented angle between $\ptl_{\omega}$ and
$\dot{\ga}(s)$.  Then we have $\dot{\omega}=\cos\sg$ and
$\dot{\tau}=(\sin\sg)/(\cos\omega)$.  Now we replace Euclidean
coordinates with spherical coordinates in the expression given above
for the mean curvature $H$ of $\Sg-R$ and we get

\begin{lemma}
\label{lem:mc4}
The generating curve $\ga=(\omega,\tau)$ in $\sph^2_+$ of a $C^2$
surface which is invariant under any rotation $r_\theta$ and has mean
curvature $H$ in $(\stres,g_h)$ satisfies the following system of
ordinary differential equations
\[
(*)_H \quad \left\{
\begin{aligned}
\dot{\omega}&=\cos\sg,
\\
\dot{\tau}\,&=\frac{\sin\sg}{\cos\omega},
\\
\dot{\sg}&=\tan\omega\,\sin\sg-\cot^3\omega\,\sin^3\sg
+2H\,\,\frac{(\sin^2\omega\,\cos^2\sg+\sin^2\sg)^{3/2}}
{\sin^2\omega},
\end{aligned}
\right.
\]
whenever $\omega\in (0,\pi/2)$.  Moreover, if $H$ is constant then the
function
\begin{equation}
\label{eq:energy}
\frac{\sin\omega\,\cos\omega\,\sin\sg}{\sqrt{\sin^2\omega\,
\cos^2\sg+\sin^2\sg}}-H\sin^2\omega
\end{equation}
is constant along any solution of $(*)_H$.
\end{lemma}

\vspace{0,1cm} Note that the system $(*)_H$ has singularities for
$\omega=0,\pi/2$.  We will show that the possible contact between a
solution $(\omega,\tau,\sg)$ and $R$ is perpendicular.  This means
that the generated surface $\Sg$ is of class $C^1$ near $R$.

The existence of a first integral for $(*)_H$ follows from Noether's
theorem \cite[\S 4 in Chap.~3]{gh} by taking into account that the
translations along the $\tau$-axis preserve the solutions of $(*)_H$.
The constant value $E$ of the function \eqref{eq:energy} will be
called the \emph{energy} of the solution $(\omega,\tau,\sg)$.  Notice
that
\[
\sin\omega\,\cos\omega\,\sin\sg=(E+H\sin^2\omega)\, \sqrt{\sin^2\omega\,\cos^2\sg+\sin^2\sg}.
\]
The equation above clearly implies
\begin{equation}
\label{eq:cuadrado}
(\sin^2\omega\,\cos^2\omega-(E+H\sin^2\omega)^2)\,
\sin^2\sg=(E+H\sin^2\omega)^2\sin^2\omega\,\cos^2\sg,
\end{equation}
from which we deduce the inequality
\begin{equation}
\label{eq:cota} \sin\omega\,\cos\omega\geq |E+H\sin^2\omega|,
\end{equation}
which is an equality if and only if $\cos\sg=0$.

Moreover, by using \eqref{eq:cuadrado} we get
\begin{equation}
\label{eq:seno}
\sin \sg = \frac{(E+H\sin^2\omega)\,\sin \omega}{\cos \omega\,
\sqrt{\sin^2\omega - (E+H\sin^2\omega)^2}}.
\end{equation}
By substituting \eqref{eq:seno} in the third equation of $(*)_H$
we deduce
\begin{equation}
\label{eq:sigmaprima} \dot{\sg}=\frac{p(\sin^2\omega)}{\cos^2
\omega\,(\sin^2\omega-(E+H\sin^2\omega)^2)^{3/2}},
\end{equation}
where $p$ is the polynomial given by $p(x)=-(E+Hx)^3-Hx^3+(E+2H)x^2$.

From the uniqueness of the solutions of $(*)_H$ for given initial
conditions we easily obtain

\begin{lemma}
\label{lem:properties}
Let $(\omega(s),\tau(s),\sg(s))$ be a solution of $(*)_H$ with energy
$E$.  Then, we have
\begin{itemize}
\item[(i)] The solution can be translated along the $\tau$-axis.
More precisely, $(\omega(s),\tau(s)+\tau_0,\sg(s))$ is a solution of
$(*)_H$ with energy $E$ for any $\tau_{0}\in\rr$.
\item[(ii)] The solution is symmetric with respect to any meridian
$\{\tau=\tau(s_0)\}$ such that $\dot{\omega}(s_0)=0$.  As a consequence,
we can continue a solution by reflecting across the critical points of $\omega(s)$.
\item[(iii)] The curve
$(\omega(s_0-s),\tau(s_0-s),\pi+\sg(s_0-s))$ is a solution of $(*)_{-H}$
with energy~$-E$.
\end{itemize}
\end{lemma}

\begin{lemma}
\label{lem:functions}
Let $(\omega(s),\tau(s),\sg(s))$ be a solution of $(*)_H$.  If
$\sin\sg(s_0)\neq 0$, then the coordinate $\omega$ is a function over
a small $\tau$-interval around $\tau(s_0)$.  Moreover
\begin{equation}
\label{eq:dx/dt}
\frac{d\omega}{d\tau}=\cos\omega\,\cot\sg,\qquad\frac{d^2\omega}
{d\tau^2}=-\frac{\sin\omega\, \cos\omega\,\sin\sg\,\cos^2\sg+\dot{\sg}\,
\cos^2\omega}{\sin^3\sg},
\end{equation}
where $\dot{\sg}$ is the derivative of $\sg$ with respect to $s$.
\end{lemma}

Now we describe the complete solutions of $(*)_H$.  They are of the
same types as the ones obtained by W.~Y.~Hsiang \cite{delaunay} when
he studied constant mean curvature surfaces of revolution in
$(\stres,g)$.

\begin{theorem}
\label{th:rotationally} Let $\ga$ be the generating curve of a $C^2$ complete, connected,
rotationally invariant surface $\Sg$ with constant mean curvature $H$ and energy $E$.  Then the
surface $\Sg$ must be of one of the following types
\begin{itemize}
\item[(i)] If $H=0$ and $E=0$ then $\gamma$ is a half-meridian and
$\Sg$ is a totally geodesic $2$-sphere in $(\stres,g)$.  \item[(ii)]
If $H=0$ and $E\neq 0$ then $\Sg$ coincides either with the minimal
Clifford torus $\mathcal{T}_{\sqrt{2}/2}$ or with a compact embedded
surface of unduloidal type.  \item[(iii)] If $H\neq 0$ and $E=0$ then
$\Sg$ is a compact surface congruent with a sphere $\mathcal{S}_H$.
\item[(iv)] If $EH\neq 0$ and $H\neq -E$ then $\Sg$ coincides either
with a non-minimal Clifford torus $\mathcal{T}_\rho$, or with an
unduloidal type surface, or with a nodoidal type surface which has
selfintersections.  Moreover, unduloids and nodoids are compact
surfaces if and only if $H/\sqrt{1+H^2}$ is a rational number.
\item[(v)] If $H=-E$ then $\ga$ consists of a union of circles
meeting at the north pole.  The generated $\Sg$ is a compact surface
if and only if $H/\sqrt{1+H^2}$ is a rational number.
\end{itemize}
\end{theorem}

\begin{figure}[h]
\centering{\includegraphics[width=10cm]{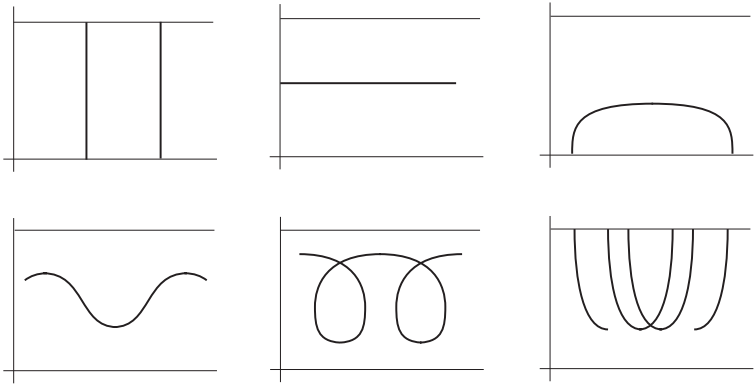}} \caption{Generating curves in spherical
coordinates of rotationally invariant surfaces with constant mean curvature in $(\stres,g_h)$.
The horizontal segments represent the $\tau$-axis and the line $\omega=\pi/2$, which is
identified with the north pole.  The vertical segment represents the $\omega$-axis.  The
generated surfaces are, respectively, a totally geodesic $2$-sphere, a Clifford torus, a
spherical surface $\mathcal{S}_\la$, an unduloidal type surface, a nodoidal type surface, and a
surface consisting of ``petals" meeting at the north pole.} \label{fig:delaunay}
\end{figure}

\begin{proof}
By removing the points where the generating curve meets the north pole
and $R$, we can suppose that $\ga=(\omega, \tau, \sg)$ is a complete
solution of $(*)_H$ with energy $E$.  By Lemma~\ref{lem:properties}
(i) we can assume that $\ga$ is defined over an open interval $I$
containing the origin, and that the initial conditions are
$(\omega_0,0,\sg_0)$.  We can also suppose that $H\geq 0$ by
Lemma~\ref{lem:properties} (iii).

To prove the theorem we distinguish several cases depending on the
value of $E$.

$\bullet$ $E=0$.  Suppose first that $H=0$.  Then $\sin \sg \equiv 0$ along $\ga$ from
\eqref{eq:seno} and so, the solution is given by $\tau\equiv 0$, $\omega(s)=s+\omega_0$ and
$\sg\equiv 0$.  We conclude that $\ga$ is a half-meridian.  The generated surface is a totally
geodesic $2$-sphere in $(\stres,g)$ with two isolated singular points.

Now suppose $H>0$.  In this case we get $\sin \sg>0$ by \eqref{eq:seno} and so we can see the
$\omega$-coordinate as a function of $\tau$.  Moreover, $\tan \omega \leq 1/H$ by
\eqref{eq:cota}, so that the solution could approach the $\tau$-axis. We can take the initial
conditions of $\ga$ as $(\arctan(H^{-1}),0,\pi/2)$.  By the symmetry of the solutions we only
have to study the function $\omega(\tau)$ for $\tau>0$.  By using \eqref{eq:sigmaprima} we
obtain $\dot{\sg}>0$, which together with the fact that $\sin \sg>0$, implies that $\sg \in
(\pi/2,\pi)$.  Therefore $\cos\sg<0$ and the function $\omega(\tau)$ is strictly decreasing. In
addition $\sin\sg\rightarrow 0$ as $\omega\rightarrow 0$ by \eqref{eq:seno} and so, $\ga$ meets
the $\tau$-axis orthogonally.

On the other hand as $\cos \sg<0$ we can see the $\tau$-coordinate as
a function of $\omega$.  This function satisfies that
\begin{displaymath}
\frac{d\tau}{d\omega}=\frac{-H\,\sin^2\omega}{\cos
\omega\,\sqrt{\cos^2\omega-H^2\sin^2\omega}},\qquad \omega
\in (0,\arctan (H^{-1})).
\end{displaymath}
We can integrate the equality above to conclude that
\begin{displaymath}
\tau(\omega)=\frac{H}{\sqrt{1+H^2}}\arcsin(\sqrt{1+H^2}\sin\omega)-
\arcsin(H\tan\omega)+\frac{\pi}{2}\,\left(1-\frac{H}{\sqrt{1+H^2}}\right).
\end{displaymath}
Finally it is easy to see that, after a translation along the $\tau$-axis, the expression of
the generated $\Sg$ in Euclidean coordinates coincides with the one given in
\eqref{eq:geocoore1} for the sphere $S_H$.

$\bullet$ $E\neq0$.  From \eqref{eq:cota} we get that
$(1+H^2)\sin^4\omega-(1-2EH)\sin^2\omega+E^2\leq 0$, which implies
that $(1-2EH)^2-4E^2(1+H^2)\geq 0$.  In this case, $\omega_1\leq
\omega\leq\omega_2$, where $\sin\omega_1$ and $\sin\omega_2$ coincide
with the positive zeroes of the polynomial
$(1+H^2)x^4-(1-2EH)x^2+E^2$.  Therefore the solution does not approach
the $\tau$-axis.  We distinguish several cases:

(i) $E>0$.  If $(1-2EH)^2-4E^2(1+H^2)=0$ ($\omega_1=\omega_2$) then
$E=(\sqrt{1+H^2}-H)/2$ and the solution is given by
$$\omega\equiv \arcsin\Big(\sqrt{\frac{1}{2}-\frac{H}{2\sqrt{1+H^2}}}\Big).$$
The generated $\Sg$ is the Clifford torus $\mathcal{T}_\rho$ with
$\rho^2=(1/2)(1+H/\sqrt{1+H^2})$.  Otherwise, by equation
$\eqref{eq:seno}$ we get that $\sin \sg >0$ and then the
$\omega$-coordinate is a function of~$\tau$.  After a translation
along the $\tau$-axis we can suppose that the initial conditions of
$\ga$ are $(\omega_1, 0,\pi/2)$.  Moreover, by symmetry of the
solutions, it is enough to study $\omega(\tau)$ for $\tau >0$.

Call $s_2$ to the first $s>0$ such that $\sg(s)=\pi/2$.  Taking into account
\eqref{eq:sigmaprima}, it is easy to see that there exists a unique $s_1\in (0,s_2)$ such that
$\dot{\sg}(s_1)=0$.  By the definition of $s_1$ and $s_2$ we get $\dot{\sg}<0$ on $(0,s_1)$ and
$\dot{\sg}>0$ on $(s_1,s_2)$, so that $\sg$ reaches a minimal value in $\sg(s_1)$.  As a
consequence, $\sg \in (0,\pi/2)$ and $\cos \sg >0$ on $(0,s_2)$.  Thus, if we define
$\tau_2=\tau(s_2)$ then the function $\omega(\tau)$ is strictly increasing on $(0, \tau_2)$ and
so, $\omega(\tau_2)=\omega_2$. On the other hand by substituting \eqref{eq:seno} and
\eqref{eq:sigmaprima} into \eqref{eq:dx/dt} we get
\begin{equation}
\label{eq:d2w/dt2} \frac{d^2\omega}{d\tau^2}=\frac{\cos\omega}
{\sin^3\omega\,(E+H\sin^2\omega)^3}\,((E+H\sin^2\omega)^3- 2(E+H)\sin^4\omega\cos^2\omega).
\end{equation}
It follows that there exists a unique value $\tau_1\in (0,\tau_2)$ such that
$(d^2\omega/d\tau^2)(\tau_1)=0$. We can conclude that the graph $\omega(\tau)$ is strictly
increasing and strictly convex on $(0,\tau_1)$ whereas it is strictly increasing and strictly
concave on $(\tau_1,\tau_2)$.  By successive reflections across the vertical lines on which
$\omega(\tau)$ reaches its critical points, we get the full solution which is periodic and
similar to a Euclidean unduloid.

As $\cos \sg >0$ on $(0,s_2)$, we can see the $\tau$-coordinate as a
function of $\omega \in (\omega_1,\omega_2)$.  Then, the period of
$\ga$ is given by
 \begin{displaymath}
T=2\int_{\omega_1}^{\omega_2}\dot{\tau}(\omega)\,d\omega=
2\int_{\omega_1}^{\omega_2}\frac{(E+H\sin^2\omega)\,\sin\omega}
{\cos\omega\,\sqrt{\cos^2\omega\,\sin^2\omega-(E+H\sin^2\omega)^2}}\, \,d\omega.
 \end{displaymath}
A straightforward computation shows that $T=(1-H/\sqrt{1+H^2})\pi$.
Then the generated $\Sg$ is an unduloidal type surfaces which is
compact if and only if $H/\sqrt{1+H^2}$ is a rational number.
Moreover, $\Sg$ is embedded if and only if $T=2\pi/k$ for some integer
$k$.  In the particular case of $H=0$, we have proved that the
generated $\Sg$ is either the minimal Clifford torus or a compact
embedded unduloidal type surface.

(ii) $E<0$.  Assuming that $H<-E$ we get that $\sin\sg <0$ along $\ga$ by \eqref{eq:seno}.  By
using the same arguments as in the previous case we deduce that $\Sg$ coincides either with the
Clifford torus $\mathcal{T}_{\rho}$ with $\rho^2=(1/2)(1-H/\sqrt{1+H^2})$, or with an
unduloidal type surface with period $T=(1+H/\sqrt{1+H^2})\,\pi$. Hence these unduloidal
surfaces are never embedded. Moreover they are compact if and only if $H/\sqrt{1+H^2}$ is a
rational number.

Thus we can suppose $H>-E$. In this case we have $\sin \sg <0$ if $\sin\omega \in
[\sin\omega_1,\sqrt{-E/H})$ while $\sin \sg
>0$ if $\sin\omega \in (\sqrt{-E/H},\sin\omega_2]$. Moreover, from \eqref{eq:sigmaprima} it is
easy to check that $\dot{\sg}> 0$ along the solution.  After a translation along the
$\tau$-axis, we can suppose that the initial conditions of $\ga$ are $(\omega_2,0,\pi/2)$.  By
the symmetry property we only have to study the solution for $s>0$.

Call $s_1$ and $s_2$ to the first positive numbers such that $s_1<s_2$, $\sg(s_1)=\pi$ and
$\sg(s_2)=3\pi/2$.  Then $\omega(s_1)=\sqrt{-E/H}$ and $\omega(s_2)=\omega_1$.  Call
$\tau_i=\tau(s_i)$, $i=1,2$.  We have that $\sg \in (\pi/2,\pi)$ on $(0,s_1)$ and $\sg \in
(\pi,3\pi/2)$ on $(s_1,s_2)$.  As a consequence, the restriction of $\ga$ to $[0,s_2]$ consists
of two graphs of the function $\omega(\tau)$ meeting at $\tau=\tau_1$. Taking into account
\eqref{eq:d2w/dt2} we can conclude that $\omega(\tau)$ is strictly decreasing and strictly
concave on $(0,\tau_1)$ whereas it is strictly increasing and strictly convex on
$(\tau_2,\tau_1)$.  As $\{\tau=0\}$ and $\{\tau=\tau_2\}$ are lines of symmetry for $\ga$, we
can reflect successively to obtain the complete solution, which is periodic.  The generating
curve is embedded if and only if $\tau_2=0$.  Let us see that this is not possible.

As $\cos\sg<0$ on $(0,s_2)$, we can see the $\tau$-coordinate as a function of $\omega$. Then,
\begin{displaymath}
\tau_2=-\int_{\omega_1}^{\omega_2}\dot{\tau}(\omega)\,d\omega=
\int_{\omega_1}^{\omega_2}\frac{(E+H\sin^2\omega)\,\sin\omega}
{\cos\omega\,\sqrt{\cos^2\omega\sin^2\omega-(E+H\sin^2\omega)^2}} \,\,d\omega.
\end{displaymath}
A straightforward computation shows that
$$\tau_2=\frac{\pi}{2}\,\Big(1-\frac{H}{\sqrt{1+H^2}}\Big)>0.$$
It follows that the period of $\ga$ is given by
$(1-H/\sqrt{1+H^2})\pi$, and $\Sg$ is a nodoidal type surface which is
compact if and only if $H/\sqrt{1+H^2}$ is a rational number.

To finish the prove we only have to study the case $H=-E$.  Now $\sin \omega_1=H/\sqrt{1+H^2}$
and $\sin \omega_2=1$.  Then the solution could approach the north pole.  Note that $\sin \sg <
0$ far away of the north pole by \eqref{eq:seno}.  Thus along any connected component of
$\ga-\{\text{pole}\}$ we can see the $\omega$-coordinate as a function of $\tau$.  Using
\eqref{eq:seno} and the expressions of $\dot{\sg}$ and $d^2 \omega/d\tau^2$ given by $(*)_H$
and \eqref{eq:d2w/dt2} respectively, it is easy to see that $\dot{\sg}> 0$ if $\omega\neq
\pi/2$ and that $d^2 \omega/d\tau^2 >0$. In addition, $\sin \sg \rightarrow 0$ as $\omega
\rightarrow \pi/2$. We can suppose that the initial conditions of $\ga$ are
$(\omega_1,0,3\pi/2)$.  By the symmetry of the solution, we only have to study $\ga(s)$ for
$s>0$.

Call $s_0$ to the first $s>0$ such that $\sin \sg (s)=0$.  As $\sg$ is strictly increasing we
get $\sg \in (3\pi/2,2 \pi)$ on $(0,s_0)$ and $\lim_{s\rightarrow s_0^-} \sg(s)=2 \pi$.  If we
call $\tau_0=\lim_{s\rightarrow s_0^-}\tau(s)<0$, we have that the function $\omega(s)$ is
strictly increasing and strictly convex on $(0,s_0)$, while $\omega(\tau)$ is strictly
decreasing and strictly convex on $(0,\tau_0)$. For $\tau=\tau_0$ the curve meets the north
pole and the tangent vector of the curve is parallel to the meridian $\{\tau=\tau_0\}$. We
continue the generating curve so that we obtain another branch of the graph of the function
$\omega(\tau)$ meeting the north pole. We can assume that $\tau(s)\to\pi+\tau_0$ modulo $2\pi$
and $\sg(s)\to\pi$ when $s\to s_0^+$. Call $s_1$ to the first $s>s_0$ such that $\sg
(s)=3\pi/2$. As $\dot{\sg}>0$ then $\sg \in (\pi,3\pi/2)$ and the function $\omega(s)$ is
strictly decreasing on $(s_0,s_1)$. Conversely, $\omega(\tau)$ is a function strictly
increasing and strictly convex on $(\tau_1,\pi+\tau_0)$ where $\tau_1=\tau(s_1)$. Note that if
$\omega(s)=\omega(\tilde{s})$ with $s\in (0,s_0)$ and $\tilde{s}\in (s_0,s_1)$, then $\sin
\sg(s)=\sin \sg(\tilde{s})$ and so, $\sg(s)+\sg(\tilde{s})=3 \pi$.  In other words, the branch
of $\omega (\tau)$ on $(\tau_1,\pi+\tau_0)$ is the reflection of $\omega(\tau)$ on $(\tau_0,0)$
across the vertical line $\{\tau=(\pi+2\tau_0)/2\}$ and $\tau_1=\pi+2\tau_0$.  By successive
reflections across the critical points of $\omega$, we obtain the full solution which is
periodic.  The solution is embedded if and only if $\tau_0=-\pi/2$. Let us see that this is not
possible.

As $\cos \sg >0$ on $(0,s_0)$ we can see the $\tau$-coordinate as a
function of $\omega$.  Then we have
\begin{eqnarray*}
\tau_0&=&\int_{\omega_1}^{\pi/2}\dot{\tau}(\omega)\, d\omega=
-\int_{\omega_1}^{\pi/2}\frac{H\sin\omega}{\sqrt{1-(1+H^2) \cos^2\omega}}\, d\omega
=-\frac{\pi}{2}\frac{H}{\sqrt{1+H^2}}>-\frac{\pi}{2}.
\end{eqnarray*}
Moreover, $\ga$ is a closed curve if and only if $\pi-2\tau_0$ is a rational multiple of
$2\pi$, which is equivalent to that $H/\sqrt{1+H^2}$ is a rational number.
\end{proof}

\begin{remark}
The surfaces described in Theorem~\ref{th:rotationally} (v) also appear in the classification
of rotationally invariant constant mean curvature surfaces in $(\stres,g)$.  However they were
not explicitly studied in \cite{delaunay}.
\end{remark}

\begin{remark}
\label{re:examplesrat}
Now we can give examples showing that all the hypotheses in
Theorem~\ref{th:irrat} are necessary.  In
Theorem~\ref{th:rotationally} we have shown that for any $H\geq 0$
there is a family of compact immersed nodoids and a family of
unduloids with constant mean curvature $H$.  As it is shown in the
proof for some values of $H$ such that $H/\sqrt{1+H^2}$ is rational
(for example $H=0$) the corresponding unduloids are compact and
embedded.
\end{remark}

\providecommand{\bysame}{\leavevmode\hbox to3em{\hrulefill}\thinspace}
\providecommand{\MR}{\relax\ifhmode\unskip\space\fi MR }
\providecommand{\MRhref}[2]{%
   \href{http://www.ams.org/mathscinet-getitem?mr=#1}{#2} }
   \providecommand{\href}[2]{#2}

\end{document}